\documentclass[11pt]{amsart}
\usepackage{graphicx}
\usepackage{amssymb}
\usepackage{amsmath}
\usepackage{amsthm}
\usepackage{epstopdf}
\usepackage[mathscr]{eucal}
\usepackage{verbatim}
\usepackage{amssymb}
\usepackage{epsfig}

\usepackage{enumerate}

\newcommand{\U}{\mathbb{U}}

\newcommand{\R}{\mathbb{R}}

\newcommand{\N}{\mathbb{N}}

\newcommand{\Q}{\mathbb{Q}}

\newcommand{\M}{\mathbb{M}}

\newtheorem{definition}{{\bf Definition}}
\newtheorem{definitions}{{\bf Definitions}}
\newtheorem{theorem}{{\bf Theorem}}
\newtheorem{proposition}{\noindent {\bf Proposition}}
\newtheorem{corollary}{\noindent {\bf Corollary}}

\newtheorem{claim}{\noindent {\bf Claim}}
\newtheorem{fact}{\noindent {\bf Fact}}

\newtheorem{remark}[definition]{\noindent {\bf Remark}}

\newtheorem{notat}[definition]{\noindent {\bf Notation}}
\newtheorem{notations}[definition]{\noindent {\bf Notations}}

\newtheorem{lemma}[definition]{\noindent {\bf Lemma}}
\def\endproof{\hfill {\kern 6pt\penalty 500
\raise -0pt\hbox{\vrule \vbox to5pt {\hrule width 5pt
\vfill\hrule}\vrule}}}

\title{Indivisible ultrametric spaces} 
\author[C.Delhomm\'e]{Christian Delhomm\'e}
\address{ E.R.M.I.T. D\'epartement de Math\'ematiques et
d'Informatique. Universit\'e de La R\'eunion, 15, avenue Ren\'e
Cassin, BP 71551,  97715 Saint-Denis Messag. Cedex 9, La R\'eunion, France} \email {delhomme@univ-reunion.fr} 
\author[C.Laflamme] {Claude Laflamme}
\thanks{The second author was supported by NSERC
of Canada Grant \# 690404}
\address{University of Calgary, Department of Mathematics and Statistics, Calgary, Alberta, Canada T2N 1N4} 
\email {laf@math.ucalgary.ca} 
\author [M.Pouzet]{Maurice Pouzet}
\thanks{This research was completed  while the third author visited the Mathematical  Department  of the University  of Calgary in summer  2006}
\address{PCS, Universit\'e
Claude-Bernard Lyon1, Domaine de Gerland -b\^at. Recherche [B], 50
avenue Tony-Garnier, 
F$69365$ Lyon cedex 07, France} \email{
pouzet@univ-lyon1.fr }
\author [N.Sauer]{Norbert Sauer}
\thanks{The fourth author was supported by NSERC of
Canada Grant \# 691325} \address{University of Calgary, Department of
Mathematics and Statistics, Calgary, Alberta, Canada T2N 1N4} \email{
nsauer@math.ucalgary.ca}

\date{\today}
\begin{document}


\begin{abstract} 
A metric space  is indivisible if for  any partition of it into finitely many pieces one piece    contains   an isometric copy of the whole space. Continuing our investigation  of indivisible metric spaces \cite{DLPS}, we show that a countable ultrametric space embeds isometrically into an indivisible ultrametric metric space if and only if  it does not contain a strictly increasing sequence of balls. \end{abstract}

\maketitle

\section*{Introduction}

A metric space $\mathbb{M}:=(M;d)$ is {\em indivisible} if for every partition   of $M$ into two parts, one of the two parts  contains an isometric copy of $\mathbb{M}$. If $\mathbb{M}$ is not indivisible then it is {\em divisible}. The notion of indivisibility was introduced for relational structures  by  R. Fra\"\i ss\'e in the fifties, see \cite{Fra} and also \cite {pouzet}, \cite {sauer}. Results obtained since then are a part of what is now called Ramsey Theory. Recently, the study of extremely amenable  groups pointed out  to indivisible metric spaces. The first step was  Pestov theorem asserting that the group $Iso(\U)$ of isometries of the Urysohn space $\mathbb U$ is extremely amenable \cite{pestov}. Next,  the discovery by Keckris, Pestov and Todorcevic \cite {KPT} of  the exact  relationship between Fraisse limits, Ramsey classes and extremely amenable groups, followed by the introduction of  the notion of {\it oscillation stable groups} and a characterization in terms of {\it $\varepsilon$-{indivisibility}}. In \cite {nesetril}, Nesetril  proving the Ramsey property of the class of ordered finite metric spaces,  suggested to look at the indivisibility properties of metric spaces. And, in \cite{hjorth}, Hjorth proved that $\mathbb U_{\Q}$, the Urysohn space with rational distances, is divisible and asked if the bounded Urysohn  $\mathbb U_{\Q \leq 1}$ is also divisible.  Prompted by the Hjorth question, we started in \cite {DLPS} to investigate indivisible metric spaces. We proved that these spaces must be  bounded and totally Cantor-disconnected (for countable spaces  a   condition  stronger than totally Cantor-disconnedness must hold,  indeed these  spaces do not contain any \emph{spider} \cite {DLPS}). This  implies that every Urysohn space $\mathbb{U}_{V}$ with a subset of $V$ dense in some initial segment of $\R_+$  is divisible, from which  the divisibility of  $\mathbb U_{\Q \leq 1}$ follows. The fact that  on every countable indivisible  metric spaces there is a natural ultrametric distance,  invited to look  at ultrametric spaces.  We proved that an  indivisible ultrametric space does not contain an infinite strictly increasing sequence of balls.  Furthermore, this condition, added to the fact that each non-terminal node in the tree associated to the space has an infinite degree,  is necessary and sufficient  for a countable homogeneous ultrametric to be indivisible \cite{DLPS}. From this follows that such a space is the ultrametric Urysohn with reversely well founded result (this latter result was also obtained by Nguyen Van Th\'e \cite {vanthe}). Here, we continue our investigation
of countable indivisible ultrametric spaces, with the idea in mind that a complete description is not out of reach. We look first at  spectra of indivisible ultrametric spaces (the spectrum of a metric space $\M:= (M, d)$ is the set  $Spec(\M):=\{ d(x,y): x,y\in M\}$).  We show that beside the fact there are subsets of $\R_{+}$ containing $0$, the only requirement imposed upon by the indivisibility is  that they have a largest element (Proposition \ref {prop:spectrum}).  Spectra of indivisible homogeneous ultrametric spaces are reversely well ordered, hence theses spaces are quite rare.  We introduce   a notion of {\it endogeneous metric  space},  generalizing the notion of homogeous metric space. We  characterize   countable endogeneous indivisible ultrametric spaces in a fashion similar to the homogeneous ones (Theorem \ref{main}). We prove that a countable ultrametric space $\M$ embeds isometrically into an indivisible ultrametric space if and only if  it does not contain an infinite  strictly increasing sequence of balls. Furthermore,  when this condition holds, $\M$ embeds into a  countable endogeneous indivisible ultrametric space with the same spectrum (Theorem \ref{thm:extindiv}). 

In Section 1 we record some  facts we will use in the rest of the paper,  the description of countable homogeneous ultrametric spaces and the special case of the indivisible ones. Except Proposition \ref {prop:spectrum}, they come from \cite {DLPS}. In Section 2 we present the notion of endogeneous ultrametric space, and  criteria  for the indivisibility of such spaces. In section 3 we present our result on the embeddability of an ultrametric space into an indivisible one.

A preliminary version of this paper was presented at the workshop on the universal Urysohn metric space, held in Beer-Sheva, May 21-24, 2006. The authors present there are pleased to thank the organizers for their warm hospitality. 

\section{Ultrametric  spaces, homogeneity and indivisibility}
We recall the following notions. 
Let $\M:=(M,d)$ be a metric space.  If $A$ is a subset of $M$, we denote by $d_{\restriction A}$ the restriction of $d$ to $A\times A$  and by $\M_{\restriction A}$ the metric space $(A, d_{\restriction A})$, that we call {\it the metric subspace of $\M$ induced on $A$}. Let  $a\in M$; for $r\in \R^+$, 
the {\it open}, resp. {\it closed},  {\it ball of center $a$, radius $r$} is the set
$B(a, r):= \{x\in M: d(a,x)<r\}$, resp. $B'(a, r):=\{x\in M:  d(a,x)\leq r\}$. For a subset $A$ of $M$, we set $B'(A, r):=\cup \{B'(a,r): a\in A\}$. In the sequel, the term {\it ball} means an open or a closed ball. When needed, we denote by $\mathcal Ball(\M)$ the collection of balls of $\M$. A ball is \emph{ non-trivial} if it has more than one element. The {\it diameter}
of a subset $B$ of $M$ is $\delta(B): =sup\{d(x,y): x,y \in B\}$. 
Four others notions will be of importance:
\begin{definitions} Let  $a\in M$, the \emph {spectrum}  of $a$ is the set $Spec(\M,a):=\{d(a, x):  x\in M\}$. The  \emph {multispectrum}  of     $\M$ is the set $MSpec(\M):= \{ Spec(\mathbb{M}, a): a\in M\}$. The \emph{spectrum} of $\M$ is the set $Spec(\M):= \bigcup MSpec(\M)$ ($=\{ d(x,y): x,y\in M\}$). The \emph{nerve} of $\M$ is the set $Nerv (\M):= \{B'(a,  r): a\in M, r\in Spec(\M,a)\}$.
\end{definitions}

\subsection{The structure of ultrametric spaces}
A metric space is an {\em ultrametric space} if it satisfies the strong triangle inequality $d(x,z)\leq \max\{d(x,y),d(y,z)\}$. See  \cite{lemin2} for example. Note that a space is an ultrametric space if and only if $d(x,y)\geq d(y,z)\geq d(x,z)$ implies $d(x,y)=d(y,z)$.
 
The essential property of ultrametric spaces is that balls are either disjoint or comparable w.r.t. inclusion. 
From this, one can look at ultrametric spaces as binary relational structures made of equivalence relations or as trees. 
\subsubsection{Equivalences  relations on ultrametric spaces}\label{basic}
 Let $\M$ be an ultrametric space. Let $x,y\in M$ and $r\in \R_{+}^{*}$, resp. $r\in \R_{+}$,  we set $x\equiv_{<r}y$, resp. $x\equiv_{\leq r}y$,  if $d(x,y)<r$, resp. $d(x,y)\leq r$. Then:
\begin{enumerate}[{(a)}]
\item The relation  $\equiv_{<r}$, resp $\equiv_{\leq r}$, is an equivalence relation; the open, resp. closed, balls of radius $r$ form a partition of $M$; the blocks of this partition being the equivalence classes of  the equivalence relation. 
\item Let  $\equiv $ be one of the equivalences  $\equiv_{<r}$, $\equiv_{\leq r}$.   Then  $x\equiv x'$, $y\equiv y'$ and $x\not \equiv y$ imply  $d(x,y)=d(x',y')$. 
\item The quotient $M/\equiv$ can be equipped with a distance $d_{\equiv}$ in such a way that  the canonical map $M\rightarrow M/\equiv$ satisfies $d_{\equiv}(p(x),p(y))=d(x,y)$ for all $x,y\in M$ such that $x\not \equiv y$.  
\end{enumerate}

\subsubsection{Valued trees}
Ultrametric spaces can be easily described in terms of   real-valued trees.  For that we recall some notions about ordered sets. Let $P$ be an ordered set (poset). We denote by $\max(P)$ the set of maximal elements of $P$. Let $x\in P$, an element $y$ of $P$ is an  {\it immediate successor}, (or a cover)  of $x$,  if $x<y$ and there is no $z\in P$ such that $x<z<y$.  0ne  usually sets $\downarrow x:=\{y\in P: y\leq x\}$ and similarly defines $\uparrow x$. We denote by ${up(P)}$ the collection of sets $\uparrow x$ where $x\in P$. The poset $P$ is a {\it forest} if $\downarrow x$ is  a chain for every $x\in P$; this is a {\it tree} if in addition every pair $x,y$ of elements of $P$ has a lower bound, and this is a  {\it meet-tree} if $x,
y\in P$ has an infimum, denoted $x\wedge y$. We say that a poset $P$ is {\it ramified} if for every $x,
y\in P$ such that $x<y$ there is some $y'\in P$ such that $x<y'$ and
$y'$ incomparable to $y$.  In the sequel, working with trees or forest, we will also use notations inherited from chains: sometimes, we will use  the notation $(\leftarrow x]$ instead of $\downarrow x$;  we will set $]a, b]:=\{x\in P: a<x\leq b\}$, $(\leftarrow a[:=\{x\in P: x<a\}$.  The poset $P$ is  {\it well founded} if every non empty subset $A$ of $P$ contains some minimal element. As it is well known, if a poset $P$ is well-founded,  for every $x, y \in P$ such that $x<y$ there is some immediate successor $z$ of $x$, such that $x<z\leq  y$. 

\begin{definition} An \emph{ultrametric tree} is a pair $(P, v)$ where $P$ is a ramified meet-tree such that every element is below some
maximal element and $v$ is a strictly decreasing map from $P$ to $\R_+$ with $v (x)= 0$ for each
maximal element $x$ of $P$. 
\end{definition}

 The following description given in \cite {DLPS} is close from the one given  by Lemin \cite{lemin2} (who
instead of $Nerv(\M)$ considered $Ball(\M)$).

\begin{theorem} \label{characterisation1}
\begin{enumerate}
\item If $\M:= (M, d)$ is an ultrametric space, then the pair $(P,
v)$, where $P: =(Nerv(\M), \supseteq)$, $\delta$) where $\delta$ is the diameter function is an ultrametric tree.
\item Conversely, if $(P, v)$ is an ultrametric tree then $\M:=(M, d)$ where $M:= \max (P)$ and $d(x,y): = v (x \wedge y)$ is an ultrametric space and $Nerv (\M)=
up(P)_{\restriction M}$ where $up(P)_{\restriction M}:=\{M\cap
\uparrow x : x\in P\}$.
\item  The two correspondences are inverse of each other. 
\end{enumerate}
\end{theorem}

In \cite{DLPS} we introduced the notion of {\it degree} of a node of a ramified meet-tree. If  $B$ is a member of the ramified meet-tree $(Nerv(\M), \supseteq)$, the degree of $B$ is the number of sons of $B$ that we define below.  
\begin{definition}
Let $\mathbb{M}:=(M;d)$ be an ultrametric space, $B\in Nerv(\M)$ and $r:= \delta(B)$.  If $r>0$,  a \emph {son} of $B$ is  any open ball of radius $r$ included into $B$; we denote by $Son(B)$ the set of sons of $B$.  \end{definition}
 Notice that according to Subsection \ref{basic}, $Son(B)$ forms a partition of $B$. Also, notice that members of $Son (B)$ do not need to belong to $Nerv(\M)$. But, if $Nerv(\M)$,  ordered by reverse of the inclusion,  is well-founded then the members of $Son(B)$ are the immediate successors of $B$ in the poset $(Nerv(B), \supseteq)$ (hence the terminology we use).

\subsection{Some examples of ultrametric spaces}\label{subsection:example}
Let $\lambda$ be a chain and let $\overline a:= (a_{\mu})_{\mu\in
\lambda }$ such that $2\leq a_{\mu} \leq \omega$. Set
$\omega^{[\overline a]}:=\{\overline b:= (b_{\mu})_{\mu\in \lambda} :
\mu\in \lambda \Rightarrow b_{\mu} <a_{\mu}$ and $supp(\overline b):=
\{\mu \in \lambda : b_{\mu} \not = 0\}$ is finite $\}$.  If
$a_{\mu}=\omega$ for every $\mu \in \lambda$, the set
$\omega^{[\overline a]}$ is usually denoted $\omega^{[\lambda]}$. Add
a largest element, denoted $\infty$ to $\lambda$. Given $\overline b,
\overline c \in \omega^{[\overline a]}$, set $\Delta(\overline b,
\overline c):=\infty$ if $\overline b= \overline c$, otherwise
$\Delta(\overline b, \overline c):=\mu$ where $\mu$ is the least
member of $\lambda$ such that $ b_{\mu}\not = c_{\mu}$.
Suppose that $\lambda$ embeds into $\R$. Let 
$w : \lambda\cup
\{\infty\}\rightarrow \R_+$ be a strictly decreasing map such that
$w(\infty)=0$, let $d_w: = w\circ \Delta $ and let $V$ be the image
of $w$. Let 
$\omega^{\leq [\overline a]}:=\{f_{\restriction (\leftarrow \mu[ }: f\in \omega^{[\overline a]},  \mu\in \lambda \cup \{\infty\}\}$ ordered by extension.   Clearly, $\omega^{\leq [\overline a]}$ is a ramified meet-tree such that every element is below some maximal element. For $\mu \in  \lambda \cup \{\infty\}$ and $f\in \omega^{\leq [\overline a]}$,  set $v(f_{\restriction\ (\leftarrow \mu[ }):= w(\mu)$. 
  
 \begin{lemma} \label{homogeneous0}\cite{DLPS}
The pair $\M:= (\omega^{[\overline a]}, d_w)$ is an ultrametric space,
$ Spec(\M)= V$ and the  ultrametric tree associated
to $\M$ is isomorphic to $(\omega^{\leq [\overline a]}, v)$.  
\end{lemma}

Let $\M$ and $\M'$ be two metric spaces.   A map $f:M\rightarrow M'$ is an {\it isometry from $\M$ into $\M'$}, or an {\it embedding},   if 
\begin{equation}\label{eqisom}
d'(f(x),f(y))=d(x,y) \; \text{for all}\; x,y\in M
\end{equation} 
This map is an isometry from $\M$ {\it onto} $\M'$ if it is surjective. We say that  $\M$ {\it embeds into} $\M'$ if there is an embedding from $\M$ into $\M'$, that $\M$ and $\M'$ {\it equimorphic} if each embeds  into the other  
 and that $\M$ and  $\M'$ are {\it isometric} if there is an isometry from $\M$ onto $\M'$. 
A \emph{local embedding from $\M$ into $\M'$} is any isometry from a subspace of $\M$ onto a subspace of $\M'$.  If $\M=\M'$, we will call it a local embedding of $\M$. 

We say that $\M$ is {\it point-homogeneous} if the group $Iso(\M)$
of surjective isometries of $\M$ acts transitively on $\M$. According to the terminology of Fra\"{\i}ss\'e \cite {Fra}, a metric space $\M$ is \emph{homogeneous} if every local embedding of $\M$ having a finite domain  extends  to an isometry of $\M$ onto $\M$ (in fact, for ultrametric spaces,  the two notions coincide\cite{DLPSHom}).

\begin{theorem}\cite{DLPS}\label {homogeneous1}
A countable ultrametric space  $\M$ is homogeneous if and only if it is isometric to some $(\omega^{[\overline a]}, d_w)$.
\end{theorem}

Let $\M$ be an ultrametric space, the {\it age} of $\M$ is the collection of finite metric spaces isometric to some subspace of $\M$. Let $V$ be a set such that $0\in V\subseteq \R_+$. Let $\mathcal
Mult _V$ (resp. $\mathcal Mult_{V, <\omega}$) be the collection of
ultrametric metric spaces (resp. finite ultrametric spaces)
$\mathbb{M}$ whose spectrum is included into $V$. 
Then $\mathcal Mult
_{V, <\omega}$ 
is closed under embeddability and has the amalgamation property.  According to the famous theorem of Fra\"\i ss\'e (1954) \cite {fraisse} p.383, if follows that if 
$V$ is countable there is a countable homogeneous ultrametric
space whose age is $\mathcal Mult _{V, <\omega}$. It has 
spectrum $V$. We denote it  $\mathbb{U}lt _V$ and we call it the {\em Urysohn ultrametric space
with spectrum $V$}.

\begin{proposition} \cite{DLPS}\label {homogeneous2}
The space $(\omega^{[\lambda]}, d_w)$ is the countable homogeneous
ultrametric space $\mathbb Ult_V$ with spectrum $V$.
\end{proposition}


\subsection{Indivisibility}
\begin{definitions}Let $\mathbb{M}:=(M;d)$ be a metric space.
The sequence $a_0,a_1,\dots, a_{n-1},a_n$ of elements of  $M$  is an {\em $\epsilon$-chain joining $a_0$ and
$a_n$} if $d(a_i,a_{i+1})\leq \epsilon$ for all $i\in n$.  Let $x,y\in M$.
 Set
\[
d^\ast(x,y):=\inf\{\epsilon>0  \text{: $x$ and $y$ are joined by an
$\epsilon$-sequence }\}.
\]
\end{definitions}

\begin{theorem}\cite {DLPS}
Let $\mathbb{M}:=(M;d)$ be a countable  homogeneous indivisible metric
space, then $\mathbb{M}^\ast$ is an homogeneous indivisible ultrametric
space.
\end{theorem}
%

\begin{theorem}\label{lem:specmax} \cite{DLPS} If an ultrametric space is indivisible then  the collection of balls,  ordered by inclusion,  is dually well-founded and the diameter is  attained. \end{theorem}

\begin{theorem}\label{thm: ultrahomind}\cite{DLPS}
Let $\M$ be a denumerable  ultrametric space.
The following properties are equivalent:
\begin{enumerate}[{(i)}]
\item $\M$ is isometric to some $\mathbb Ult_V$, where $V$ is dually  well-ordered.
\item $\M$ is point-homogeneous, $Nerv(\M)$ ordered by reverse of the inclusion is well founded and every non-trivial $B\in Nerv(\M)$ has infinitely many sons. 
\item $\M$ is  homogeneous and indivisible. 
\end{enumerate}
\end{theorem}

This result (in part)  was obtained independently by L. Nguyen Van Th\'e \cite{vanthe}.

The crucial part is the implication $(ii)\Rightarrow (iii)$. It  is now a consequence of Theorem \ref{main}.
\subsection{Spectrum of indivisible ultrametric spaces}
The proposition below could be derived from Theorem \ref{thm:extindiv}. The proof we give uses only   Theorem  \ref {thm: ultrahomind}.
\begin{proposition}\label{prop:spectrum}
A set  $V$ is the spectrum of an ultrametric space if and only if $0\in V\subseteq \R_{+}$. If this latter condition is fulfilled, $V$ is the spectrum of an indivisible ultrametric space if and only if  $V$ has a largest element. In this case $V$ is the spectrum of an indivisible ultrametric space of size $\aleph_0+\vert V\vert $.
\end{proposition}

\begin{proof}
If $V=Spec(\M)$ for some metric space $\M$ then clearly $0\in V\subseteq\R_{+}$. Conversely, let $V$ such that $0\in V\subseteq \R_{+}$. Define $d:V\times V\rightarrow \R_{+}$, setting $d(x,y):= \max \{x, y\}$ if $x\not =y$ and $d(x,y):=0$ otherwise. Then $\M:= (V,d)$ is an ultrametric space for which $Spec(\M)=V$.  If $\M$ is an indivisible ultrametric space, its diameter is attained (Theorem \ref{lem:specmax}), that is $V:= Spec (\M)$ has a largest element. Conversely, let $V$,  with a largest element $r$,  such that  $0\in V\subseteq \R_{+}$.  We set $\M:= \U lt_V$,  the Urysohn ultrametric space with age  $\mathcal M_{V,<\omega}$,  if $V$ is finite. Otherwise, let $M:= \bigcup \{Ult_F\times \{F\}:F\in D\}$  where $D$ is the set of finite  subsets $F$ of $V\setminus \{r\}$ containing $0$. For two elements $(x,F), (x',F')\in M$ set $d((x, F),(x',F')):= r$ if $F\not =F'$, otherwise set $d((x, F),(x',F')):= d(x,x')$ where $d$ is the distance on $\U lt_F$. Clearly, $\M$ is an ultrametric space with spectrum $V$.  If $V$ is finite,  $\M$ is indivisible by Theorem  \ref {thm: ultrahomind}. 
Suppose that $V$ is infinite. Let $f: M\rightarrow 2$. Set $g(F)=0$ if there is some isometry $\psi_{F, 0}: \mathbb U lt_F\rightarrow \M_{\restriction Ult_F\times \{F\}\cap f^{-1}(0)}$. Otherwise, set $g(F)=1$, and since by Theorem 
\ref {thm: ultrahomind},  $\mathbb U lt$ is indivisible, select   an isometry $\psi_{F, 1}: \mathbb U lt_F\rightarrow \M_{\restriction Ult_F\times \{F\}\cap f^{-1}(1)}$.  Ordered by inclusion,  $D$ is up directed. It follows that for some $i<2$, $g^{-1}(i)$ is  {\it cofinal} in $D$, that is every member of $D$ is included into some member of $g^{-1}(i)$. In fact, as it is easy to see, more is true: there is a one to one map $\varphi : D\rightarrow g^{-1}(i)$ such that $F\subseteq \varphi (F)$ for every $F$ in $D$. Since $\U lt_F $ embeds into $\U lt _{F'}$ by some map $e_{F,F'}$ whenever $F\subseteq F'$, we may define a  map $\psi: M\rightarrow M$  by $\psi((x,F)):=(\psi_{\varphi(F), i}(e_{F,\varphi(F)}(x)), \varphi(F))$. This map  is an isometry from $\M$ into $\M_{\restriction f^{-1}(i)}$ proving that $\M$ is indivisible.
\end{proof}

\section{Endogeneous and indivisible ultrametric spaces}

\subsection{Endogeneity}
 \begin{definition} Let $\M$ and $
\M'$ be two metric spaces, a  \emph{local spectral-embedding}, in brief a \emph{local spec-embedding}, is a local embedding  from $\M$ into $\M'$ such that:
\begin{equation} Spec(\M, x)\subseteq Spec(\M', f(x)) \; \text{for every}æ\;  x\in Dom(f).  
\end{equation}
If $\M=\M'$ we will simply speak of  local spec-embedding   \emph{ of}  $\M$. 
\end{definition}

\begin{definitions}
Let $\M$ be a  metric space.
\begin{enumerate}[{(a)}]
\item  $\M$ is  \emph{spec-endogeneous} if every   local spec-embedding  of $\M$ extends to an embedding of $\M$.
\item $\M$ satisfies the \emph {spec-extension property}  if  for every $y\in M$, every local spec-embedding $g$ of $\M$ defined on $y$ extends to every other element $x$ to a local spec-embedding of $\M$. 
\item If furthermore, there are infinitely many  such extensions to $x$ whose images are pairwise at distance $d(x,y)$, then $\M$ satisfies the \emph{infinite extension property}.
\end{enumerate}
\end{definitions}

 \begin{notations} Let $x\in M$, $r\in \R_{+}$ and  $B$ be a ball. We set $S(x, r):= \{y\in M: d(x, y):=r\}$, $M(x):= \{y\in M: Spec(\M, x)\subseteq Spec(\M, y)\}$, $M(\neg x):= M\setminus M(x)$ and  $B(x):= B\cap  \M(x)$
\end{notations}
With these notations, the  definition (c) above requires that  $M(x)\cap S(g(y), d(x,y))$ contains an infinite set whose elements are pairwise at distance $d(x,y)$. 

\begin{lemma}\label{lem:smalltrick}
If $x\in B$ then $B(x)= \{y\in B: Spec(\M_{\restriction B}, x)\subseteq  Spec(\M_{\restriction B}, y)\}$
\end{lemma}
We have easily:

\begin{lemma}\label {lem:break}  Let $\M$ be an ultrametric space. The following properties are equivalent:
\begin{enumerate}[{(i)}]
\item $\M$  satisfies the infinite spec-extension property.\label{infiniteext}
\item \begin{enumerate}
       \item  $\M$ satisfies the spec-extension property.
       \item  For every $x, y\in M$, with $x\not =y$, the set $C_{y, x}:=M(x)\cap S(y, d(y,x))$ contains infinitely many elements at distance $d(x,y)$ from each other. 
\end{enumerate}
\end{enumerate}
\end{lemma}

\begin{proposition} \label{thmendendogeneous1} A  countable metric space $\M$ satisfying the infinite extension property is spec-endogeneous.
\end{proposition}
\begin{proof}
We prove by induction on $n$ that  every local spec-embedding $f$ of  $\M$,   with domain $A$ having  size   at most $n$,  extends to every $x\in M\setminus A$ to a local spec-embedding $\overline f$ of  $\M$. Since $M$ is countable and every increasing union of local spec-embedding is a spec-embedding, this will insure  that $\M$ is spec-endogeneous. 
Let $n<\omega$, $A\subseteq M$, with $\vert A\vert =n$  and $x\in M\setminus A$. If $n=0$, the identity map provides the required extension. Suppose $n>0$. Set $r:= d(x,A):= \min \{d(x,y):= y\in A\}$   and $A_0:=\{y\in A: d(x,y)=r\}$. Let  $y\in A_0$. Since $d(x,y)=r$, $r\in Spec(\M,y)$ and, since $f$ is a local spec-embedding, $Spec(\M, y)\subseteq Spec(\M, f(y))$, hence $r\in spec(\M, f(y))$, that is $B':=  B'(f(y),r)\in Nerv(\M)$. Since $f$ is an isometry on $A$,   $B'$ is independent of $y$. 
%

Pick  $y_0\in A_0$. Since $\M$ satisfies the infinite spec-extension property,  the set $C:= M(x)\cap S(f(y_0), r))$ contains infinitely many elements pairwise a distance $r$.  The set $\bigcup \{B(f(y),r):y\in A_0\}$ contains no more than $\vert A_0\vert $ elements at distance $r$, hence it does not cover $C$. Pick $x' \in C\setminus \bigcup \{B(f(y),r): y\in A_0\}$. Extend $f$ by setting $\overline{f}(x):=x'$. 

\begin{claim}\label{claimspec}
$\overline{f}$ is a spec-embedding.
\end{claim} 
\noindent {\bf Proof of Claim \ref{claimspec}. }
This claim amounts to:
\begin{enumerate} 
\item $d(\overline {f}(x), f(y))=d(x,y)$ for all $y\in A$.
\item $Spec (\M, x)\subseteq Spec(\M, \overline {f}(x))$.
\end{enumerate}
Item (1). Let $y\in A$. Set $r':= d(x,y)$. If $y\in A_0$, $r'=r$. Since $\overline {f}(x)\not \in B(f(y),r)$ and $C\subseteq B'$,  $d(\overline {f}(x), f(y))=r'$. If $y\in A\setminus A_0$, then by the definition of $r$, $r'>r$. Since $d(x,y_0)=r$, we have $d(y_0,y)=r'$ hence $d(f(y_0), f(y))=r'$. Since $d(\overline {f}(x),f(y_0))=r$, it follows that $d(\overline {f}(x), f(y))=r'$, as required. 

\noindent Item (2). This follows from the fact that $x'\in C$. 
\endproof


With Claim \ref{claimspec} the proof of  Proposition \ref {thmendendogeneous1} is complete. 
\end{proof}
\begin{corollary}\label{prop:sons}
For a countable  ultrametric space $\M$ the following properties are equivalent:
\begin{enumerate}[{(i)}]
\item $\M$ satisfies the infinite spec-extension property.
\item \begin{enumerate}
\item  $\M$ satisfies the spec-extension property.
\item For every $B\in Nerv(M)$ and every son $B'$ of $B$ there are infinitely many sons $B''$ such that $\M_{\restriction B'}$ embeds  into $\M_{\restriction B''}$.
\end{enumerate}
\end{enumerate}\end{corollary}
\begin{proof}
If Property (ii)b holds then Property  (ii)b of Lemma \ref{lem:break}  holds. Indeed, let $x, y\in M$. Set $r:= d(x,y)$ and   $B:= B'(y, r)$ Then $B\in Nerv(\M)$ and $B':= B(x, r)$ is a son of $B$. Moreover,  if $x'$ is  the image of $B'$ by some embedding $g$ into $B$ then $Spec(\M, x) \subseteq Spec(\M, x')$ (Lemma \ref {lem:smalltrick}).
Now, if 
 $x',x''$ are the images of $x$ into two distinct sons, then $d(x',x'')=r$. Hence $C_{y, x}$ contains infinitely many elements, as claimed.  Thus with   the spec-extension property the infinite spec-extension holds. For the converse, let $B\in Nerv(\M)$ and $B'$ be a son of $B$. Pick $x\in B'$, $y\in B\setminus B'$. Property  (2)b of Lemma  \ref{lem:break} 
asserts that $x$ can be spec-embedded into infinitely many sons of $B$. Since $\M$ is countable, Proposition \ref{thmendendogeneous1} applies and $M_{\restriction B'}$ embeds into these sons. 
\end{proof}

\subsection{Multispectrum, endogeneity and indivisibility} 

\begin{proposition}\label{prop:updirected}
Let $\M$ be a countable metric space such that every non trivial member of $Nerv(\M)$ has infinitely many sons. Then the following properties are equivalent:

\begin{enumerate} [{(i)}]
\item
\begin{enumerate}
\item Every local spec-embedding of $\M$ defined on a singleton extends to an embedding of $\M$.
\item For every $B\in Nerv(\M)$, $MSpec (\M_{\restriction B})$ is up-directed.
\end{enumerate}

\item 
\begin{enumerate} 
\item For every $y, y', x\in M$, if  $Spec(\mathbb{M}, y)\subseteq Spec(\mathbb{M}, y')$, there is some $x'\in B'(y', d(x,y))$ such that $Spec( \M, x)\subseteq Spec(\M,x')$. 
\item For every  $B\in Nerv(\M)$ and  every $a\in B$, $\M_{\restriction B}$ embeds  into $\M_{\restriction B(a)}$. 
 \end{enumerate}
 \end{enumerate}
%

%

%


\end{proposition}
\begin{proof}
Suppose that (i) holds.  
\noindent First  (ii)(a) holds trivially. 
\noindent Next (ii)(b) holds. 
For that we prove first that $\M$ has the infinite spec-extension property. We use Lemma \ref{lem:break}. Let $x, y\in M$, with $x\not =y$. Let $\mathcal B_{x}:= \{B'\in Son(B): B'\cap M(x)\not = \emptyset\}$ and $C_{y, x}:=M(x)\cap S(y, d(y,x))$. Clearly $C_{y, x}$ contains infinitely many elements at distance $d(x,y)$ from each other if and only if $\mathcal B_{x}$ is infinite. Suppose that $\mathcal B_{x}$ is finite. Let $B'\in Son(B)\setminus \mathcal B_{x}$. Pick $x'\in B'$. Since $MSpec (\M_{\restriction B})$ is up-directed, there is some $z\in B$ such that  $Spec(\M, x)\cup Spec(\M, x')\subseteq Spec(\M, z)$ (use Lemma \ref{lem:smalltrick}).  Since (i)(a) holds,  there is an embedding $f$ of $\M$ such that $f(x')=z$. This embedding maps each member of $\mathcal B_{x}$ into a member of $ \mathcal B_x$, and $B'$ into a member of $ \mathcal B_x$. This contradict the supposed finiteness of $\mathcal B_x$. 

Next, let $a\in B$. We prove by induction on $n$ that  every local spec-embedding $f$ of  $\M_{\restriction B}$,   with domain $A$ having  size   at most $n$  and range  included into $B(a)$,  extends to every $x\in B\setminus A$ to a local spec-embedding $\overline f$ of  $\M_{\restriction B}$ with $\overline{f}(x)\in B(a)$. Since $B$ is countable, this will insures that $\M_{\restriction B}$ embeds into $\M_{\restriction  B(a)}$.  We do exactly as is the proof of Proposition  \ref{thmendendogeneous1}. At the final stage, we only have to check that  the set $D:= B(x)\cap B(a) \setminus  \cup \{B(f(y), r): y\in A_0\}$ is non empty. Since  $MSpec (\M_{\restriction B})$ is up-directed, there is some $c\in B$ such that $B(c)\subseteq B(x)\cap B(a)$. Since, from  the proof of $ii(a)$ above, $\mathcal B_c$ is infinite, $D$ is nonempty.  
%



%
Conversely, that (ii) holds. 
(i)(b)  follows easily (i)(b). To get that (i)(a) holds it suffices from suppose that $\M$ satisfies properties $(a)$, $(b)$ and $(c)$.  From $(b)$ $MSpec (\M_{\restriction B})$ is up-directed, that is  property (2) holds. To conclude, it suffices to prove that $\M$ has the  infinite spec-extension property and to apply Proposition \ref{thmendendogeneous1}. 

For that, let  $y, y', x\in M$ such that  $Spec(\mathbb{M}, y)\subseteq Spec(\mathbb{M}, y')$. Set $r:= d(x,y)$, $B':= B'(y', r)$, $C':= \{x'\in B': Spec (\M, x)\subseteq Spec(\M, x')\}$. Our aim is to show that $C'\cap \{ x'\in B': d(y', x')=r\}$ contains infinitely many elements at distance $r$  of each other. This amounts to show that $C'$ has this property.  From $(c)$, $C'$ is non empty. Let $a\in C'$.  From $(b)$,   $\M_{\restriction B'}$ embeds into $\M_{\restriction B'(a)}$. According to $(a)$, $B'$ contains infinitely many elements at distance $r$ of each other.  Since  $B'(a)\subseteq C'$, $C'$ enjoy this property too. 
\end{proof}


\begin{lemma} \label{updirectedindiv}If $\M$ is indivisible then 
\begin{enumerate}
 \item $M \in Nerv(\M)$ and for each son $B$ of $M$ there are infinitely many sons $B'$ such that $\M_{\restriction B}$ embeds  into $\M_{\restriction B'}$.
\item For every $x\in M$, $\M$ embeds into $\M_{\restriction M(x)}$.
\end{enumerate}
\end{lemma} 
\begin{proof}
\noindent Item (1). The fact that $M\in Nerv(\M)$ follows from Theorem \ref {lem:specmax}. Let $r:= \delta(\M)$. If $r=0$, $M$ has no son and the property holds. So we may suppose $r\not=0$. Since $M\in Nerv(\M)$, $r$ is attained, hence $M$ has at least two sons.  Let $B\in Son(M)$. Suppose that $M$  has only finitely many sons $B_1, \dots, B_k$ such that $M_{\restriction B}$ embeds into $\M_{\restriction B_i}$ for $i=1, \dots , k$. Let $\mathcal B:= \{B'\in Son(M): \M_{\restriction B} \; \text{does not embed into}   \;  M_{\restriction B'}\}$ and $B_0:= \cup \mathcal B$.    The sets $B_0, \dots, B_k$ form a partition of $M$. Since  $M$ is indivisible, $\M$ embeds into some $\M_{\restriction B_i}$. Since $\delta(\M) >\delta(M_{\restriction, B_i})$ for $i>0$, we have $i=0$. But this  is impossible, indeed, if $g$ was an embedding, it would send two elements $x$ and $y$ of $B$  into two different sons and we would have $d(x,y)<r=d(g(x), g(y))$.  
%

\noindent Item (2). We have $M=M(x)\cup M(\neg x)$. Trivially, $\M$ does not embeds into $\M_{\restriction M(\neg x)}$. The conclusion follows with  the indivisibility of $\M$.  
\end{proof}

\begin{definition}
A metric space $\M$ is \emph {hereditarily indivisible}  if $\M$ is indivisible and for every ball $B$, $\M_{\restriction B}$ is indivisible.  \end{definition}

We get for spec-endogeneous metric spaces the  analog of the equivalence $(ii)\Leftrightarrow (iii)$ of Theorem \ref {thm: ultrahomind}. 
\begin{theorem}\label{main} A  countable ultrametric space $\M$ is spec-endogeneous and hereditarily indivisible if and only if it satisfies the following properties:
\begin{enumerate}

\item Every local spec-embedding of $\M$ defined on a singleton extends to an embedding of $\M$.
\item $(Nerv(\M), \supseteq)$ is well founded. 
\item Every non-trivial ball of $Nerv(\M)$ has infinitely many sons.
\item 
For every  ball $B\in Nerv(\M)$, $MSpec (\M_{\restriction B})$ is up-directed. 
 \end{enumerate}
\end{theorem}
%

%

\begin{proof}
Suppose that $\M$ is spec-endogenous and hereditarily indivisible. We prove successively that properties (1) , (2), (3)  and  (4) are satisfied. 

\noindent Item (1). Follows from the fact that $\M$ is spec-endogeneous.

\noindent Item (2). Follows from the fact that $\M$ is indivisible,  with the help of Theorem \ref {lem:specmax}.

\noindent Item (3). Since $\M$ is spec-endogeneous,  it has the spec-extension property. Since it is hereditarily indivisible, each non-trivial ball in $Nerv(\M)$ embeds into infinitely many sons  (Lemma \ref {updirectedindiv}).

\noindent Item (4). Follows from the fact that $ \M_{\restriction B}$ is indivisible with the help of Lemma \ref{updirectedindiv}.

Conversely, suppose that $\M$ satisfies properties $(1)$, $(2)$, $(3)$, $(4)$. 
First, from $(1)$, $(3)$ and $(4)$, $\M$ has the infinite spec-extension property (Proposition \ref{prop:updirected}). Since $\M$ is countable,  $\M$ is spec-endogeneous (Proposition \ref{thmendendogeneous1}). To conclude, we have to  show that  $\M$ is hereditarily indivisible.  It suffices to prove that $\M$ is indivisible.  Indeed, if $B\in Nerv(\M)$, $\M_{\restriction B}$ satisfies $(1)$, $(2)$, $(3)$ and $(4)$.  Hence, by the same token,  $\M_{\restriction B}$ will be indivisible. 

\begin{claim}\label{claim:cofin}
For each non-trivial  $B\in Nerv(\M)$ and every finite set $\mathcal C$ of sons of $B$,
$\M_{\restriction B}$ embeds into $\M_{\restriction B\setminus \cup\mathcal C}$.
\end{claim}

\noindent {\bf Proof of Claim \ref{claim:cofin}.}
Since $\M$ has the infinite  spec-extension property, for every ball $B$ of $\M$, $\M_{\restriction B}$ has this property.  Let $B\in Nerv(\M)$. From Corollary \ref{prop:sons},  for  every $B'\in Son(B)$ there are infinitely many $B''\in Son(B)$ such that $\M_{\restriction B'}$ embeds into $\M_{\restriction B''}$. Since   $Son(B)$ is countable, there is a one-to-one mapping $\psi: Son(B)\rightarrow  Son(B) \setminus\mathcal C$ such that  $\M_{\restriction B'} $embeds into  $\M_{\restriction \psi(B')} $ for each $B\in Son (B)$.  With the fact that $B= \cup Son (B)$, this implies that $\M_{\restriction B}$ embeds into  $\M_{\restriction B\setminus \mathcal C}$.
\endproof

Let  $\chi: M\rightarrow\{0,1\}$ be a bicoloring of $M$. 
Let $\mathcal M_0$ denote the set of balls $B\in Nerv(\M)$ such that there is some
isometry $\varphi_B$ from $\M_{\restriction B}$ into  $\M_{\restriction B\cap\{\chi^{-1}(0)\}}$ and let $M_0:=\cup\mathcal M_0$.
Observe that $M_0\supseteq \chi^{-1}(0)$. 

\begin{claim}\label{claim-retrM0}

\begin{enumerate}
  \item For every subset $\mathcal N$ of $\mathcal M_0$, there is an isometry of $\M_{\restriction \cup
  \mathcal N}$ into $\M_{\restriction (\cup\mathcal N)\cap\chi^{-1}(0)}$.
  \item  Let $B\in Nerv(\M)$.
 If $B$ is included in no member of $\mathcal M_0$, then $\M_{\restriction B}$  does not embed in $\M_{\restriction B\cap M_0}$.\end{enumerate}
\end{claim}

\noindent {\bf Proof of Claim \ref {claim-retrM0}.}
Both parts rely on the fact that balls are either disjoint or comparable w.r.t. inclusion. 
\begin{enumerate}
  \item 
  Let $\mathcal N'$ denote the set of maximal members of $\mathcal N$ (maximal w.r.t. inclusion). 
Let $\varphi:=\cup \{\varphi _B: B\in\mathcal N'\}$.  Since balls are either disjoint or comparable, 
$\varphi$ is a function and, since $P:= (Nerv(\M), \supseteq)$ is  well-founded,
$\cup\mathcal N'= \cup\mathcal N$, hence the domain of $\varphi$ is $\cup\mathcal N$.
   \item Since $B$ is assumed to be included in no member of $\mathcal M_0$, and balls are either disjoint or comparable, 
  $B\cap M_0=B\cap(\cup\mathcal M_0)=\cup\mathcal\{B'\in \mathcal M_0: B'\subseteq  B\}$. 
Hence, according to the first part of the present claim, $\M_{\restriction B\cap M_0}$ embeds into  $\M_{\restriction B\cap \chi^{-1}(0)}$. On the other hand $\M_{\restriction B}$ does not embed into $\M_{\restriction B\cap \chi^{-1}(0)}$,
since we have supposed that  $B\notin\mathcal M_0$. It follows that $\M_{\restriction B}$ does not embed into $\M_{\restriction B\cap M_0}$.
\end{enumerate}
\endproof
 
Now suppose  that $M\notin\mathcal M_0$.  

\begin{claim}\label{claim:main} Every local spec-embedding $f$ of  $\M$ with  a finite domain $A$ and  its range included into $M\setminus M_0$ extends to every $x\in M\setminus A$ to a local spec-embedding $\overline f$ of    
$\M$ with range included into $M\setminus M_0$.

 \end{claim}
\noindent{\bf Proof of Claim \ref{claim:main}.} We argue by induction on $n:= \vert A\vert $. We proceed as  for the proof of Proposition \ref{thmendendogeneous1}.  Suppose $n=0$. Since $Nerv(\M)$ is well-founded, $M\in Nerv(\M)$ and we may apply  Proposition \ref{prop:updirected}. Thus  $\M$ embeds into $\M_{\restriction M(x) }$. Since $\M$ does not embed into  $\M_{\restriction {M_0}}$ (Claim~\ref{claim-retrM0}), $M(x)\setminus  M_0$ is non empty~; choose  any element $x'$ in it and set $\overline f(x):= x'$.  

Suppose $n>0$. Set $r:= d(x,A):= \min \{d(x,y):= y\in A\}$, $A_0:=\{y\in A: d(x,y)=r\}$ and $\mathcal C:= \{B(f(y), r): y\in A_0\}$.  Our aim is to find some $x'$ in the intersection of $M\setminus M_0$, $M(x)$ and $\cap \{S(f(y), r): y\in A_0\}$. Indeed, setting $\overline f(x):=x'$, the same argument as in Proposition \ref{thmendendogeneous1} yields that $\overline f$ is a spec-embedding.

Let  $y\in A_0$. Since $d(x,y)=r$, $r\in Spec(\M,y)$ and, since $f$ is a local spec-embedding, $Spec(\M, y)\subseteq Spec(\M, f(y))$, hence $r\in spec(\M, f(y))$, that is $B':=  B'(f(y),r)\in Nerv(\M)$. Since $f_{\restriction A_{0}}$ is an isometry, $B'$ is independent of $y$. 

Our aim reduces to find some $C\in Son(B')\setminus \mathcal C$ and such that $C(x):= C\cap M(x)$ is not included into $M_0$.  For that, it suffices to prove that 
\begin{equation} \M_{\restriction B'}\;   \text {embeds into}\;   \M_{\restriction B'(x)\setminus \cup \mathcal C}\end{equation}
 Indeed, according to Claim~\ref{claim-retrM0}, $\M_{\restriction B'}$ does not embed into $\M_{\restriction B'\cap N_0}$. Hence $B'(x)\setminus \cup \mathcal C$ is not included into $M_0$. Since $B'(x)\setminus \cup \mathcal C=\cup \{C(x): C\in Son(B)\setminus \mathcal C\}$, there is some $C\in Son(B)\setminus \mathcal C$ such that $C(x)$ is not included into $M_0$. 

To get $(4)$, we prove  first that  $\M_{\restriction B'}$æembeds into $\M_{\restriction B'(x)}$. Indeed, 
pick  $y_0\in A_0$. Since $\M$ satisfies the infinite spec-extension property, there is some $x'\in M$ such that $d(x',f(y_0))=d(x, y_0)$ and $Spec(x, \M)\subseteq Spec(\M, x')$.  Let $B'':= \{x'\in B' : Spec(\M_{\restriction B'}, x )\subseteq Spec(\M_{\restriction B'}, x')\}$. We have $B''=B'(x)$, hence, from  Proposition \ref {prop:updirected}, $\M_{\restriction B'}$æembeds into $\M_{\restriction B'(x)}$. Next, applying  Claim \ref{claim:cofin} we get that 
 $\M_{\restriction B'}$æembeds into $\M_{\restriction B\setminus \cup \mathcal C}$. If  $g$ and $ h$ are  two such embeddings, in this order, then $h\circ g$ is an embedding of $ \M_{\restriction B'}$ into $\M_{\restriction B'(x)\setminus \cup \mathcal C}$. 
\endproof

Since $M$ is countable, Claim \ref{claim-retrM0} insures that $\M$ embeds into  $\M_{\restriction M\setminus M_0}$.  Since $M\setminus M_0\subseteq \chi^{-1}(1)$, $\M$ embeds into $\M_{\restriction \chi^{-1}(0)}$. This proves that $\M$ is indivisible.

\end{proof}

\section{Extensions of indivisible ultrametric spaces}
The purpose of this section is to prove:

\begin{theorem}\label{thm:extindiv}
 A countable ultrametric space $\M$ embeds into  a countable  indivisible ultrametric space if and only if it does not contain an infinite strictly increasing sequence of balls. Furthermore, when this condition holds $\M$ embeds into a countable  spec-endogeneous indivisible ultrametric space with the same  spectrum \end{theorem}

The fact that the condition on balls is necessary follows from  Theorem \ref{lem:specmax}. For the sufficiency, we construct an extension of $\M$  to which we  can apply Theorem \ref{main}. 

The key notions  are  these: 
\begin{definitions} Let $\M$ be a metric space; a binary operation, denoted $+$, on $\M$ is \emph{ compatible}
if 
\begin{equation} \label{eqcompatible1} d(z+x, z+y)= d(x, y) =d(x+z, y+z)
\end{equation} for all $x,y,z \in M$. 
An \emph{ultrametric monoid} is an ultrametric space $\M$ endowed  with a compatible operation $+$ such that $M$ with this  operation is a monoid.
\end{definitions}
Indeed, we will prove that $\M$ extends to an ultrametric monoid $\M^*$ with the same spectrum, and the same condition on balls (Theorem \ref{thm:firstext}) and having infinitely many sons. Next, we will extend $\M^*$ to an other ultrametric monoid, $Path(\M^*)$, such that each of its  balls can be also endowed with a structure of ultrametric monoid (Theorem \ref{thm:pathext}, and Theorem \ref{thm:compat}). Finally we will prove that this space is spec-endogeneous and hereditary indivisible (Theorem \ref {doubleextension1}).

\subsection{Monoids extensions of an ultrametric space} 

Let $\M:= (\omega^{[\lambda]}, d_w)$, as defined in Subsection\ref{subsection:example}. For $f,g\in \omega^{[\lambda]}$ let $f+g$ be defined by $(f+g)(\mu)=f(\mu)+g(\mu)$, and let $0$ be the constant map equal to $0$. With this, $M$ is a comutative  monoid. Furthermore, the operation is compatible. Thus $\M$ is a commutative ultrametric monoid.

The set   $\omega ^{\leq [\lambda]}$  ordered by extension is a ramified mee-tree in which every element is below some maximal one. For $x,x'\in \omega ^{\leq [\lambda]}$, we denote by $x\wedge x'$ the meet of  $x,x'$. Let $X\subseteq \omega^{[\lambda ]}$. Set $X^*$ for the set of finite sums of members of $X$ with $0$ included. Let $T(X):= \{e\wedge e': e,e'\in X\}$.  It is easy to show that  $T(X)$ is a meet-tree; we call it \emph{the meet-tree generated by $X$}. 
\begin{lemma}
Let $X\subseteq \omega^{[\lambda ]}$.  If $T(X)$ is well-founded then $T(X^{*})$ is well-founded.
\end{lemma}
\begin{proof}
Suppose that $T(X^*)$ is not well-founded. Let $Y\subseteq X^*$ and let 
$y_0>y_1\cdots>y_n>\cdots$  be  an infinite strictly decreasing sequence of members of $T(Y)$. Let $\overline y\in Y$ such that $\overline y\geq y_0$.  

{\bf Claim 1.} There is an infinite sequence $\overline y_0,\overline y_1,\dots \overline y_n, \dots$ of members of $Y$ such that $y_n= \overline y\wedge \overline y_n$ for all $n\in \N$.

{\bf Proof of Claim 1. }   Since $y_n\in T(Y)$ there are $e_n, e'_n \in Y$ such that $y_n=e_n\wedge e'_n$. Since  $\overline y\geq y_n$, we have either $y_n=\overline y\wedge e_n$ or $y_n=\overline y\wedge e'_n$. In the first case  set $\overline y_n:=e_n$ and in the second case $\overline y_n:=e'_n$.  \endproof

With no loss of generality, we may suppose $\overline y> y_0$ (otherwise a subsequence of the $y_n$'s will do). Thus for every $n<\omega$,  $\overline y\not =\overline y_n$. Let $\lambda_n$ be the least  element of $\lambda$ such that 
\begin{equation}\label {eqtree1}
\overline y(\lambda_n)\not = \overline y_n(\lambda_n)
\end{equation}
 Then $(\leftarrow \lambda_n[$ is the domain of $\overline {y}\wedge \overline {y}_n$. 

Since $y_0>y_1>\cdots>y_n>\cdots$ we have $\lambda_0>\lambda_1>\cdots>\lambda_n>\cdots$. 

Let $A:= \cup\{]\lambda_{n+1} , \lambda _n[ :n<\omega\}$, $B:= \cap  \{(\leftarrow , \lambda _n[ :n<\omega\}$,  

Since $Supp (\overline y)$ is finite, we may suppose that  it is  disjoint from $A$ (otherwise a subsequence of the $\overline y_n$'s will do). In particular $\overline y(\lambda_n)=0$.  From (\ref{eqtree1}) we have $\overline y(\lambda_n)\not = \overline y_n(\lambda_n)$. Hence $\overline y_n(\lambda_n)>0$.

Since $\overline y_n\in Y\subseteq X^*$, this is a finite sum  of members of $X^*$. We may choose such a member $\overline z_n$ such that  $\overline z_n(\lambda_n)\not =0$. Since $\overline z_n$ is a term  in a finite sum which is equal to $\overline y_n$, we have $\overline z_n\leq \overline y_n$, meaning that 

\begin{equation}\label{ineqtree2}
\overline z_n(\mu)\leq \overline y_n(\mu) \; \text {for every} \;  \mu<\lambda. 
\end{equation}

 {\bf Claim 2.} $\overline z_n$ is $0$ on $(\leftarrow \lambda_n[\cap A$. 

{\bf Proof of Claim 2.} 

Combine (\ref{eqtree1}) and  inequality (\ref {ineqtree2}) for $\mu \in (\leftarrow \lambda_n[\cap A$.
\endproof

{\bf Claim 3.}  There are infinitely many $\overline z_n$'s whose restrictions to $B$ are the same.

{\bf Proof of Claim 3.} 
Since $B\subseteq (\leftarrow \lambda_n[)$, $\overline y_n$ and $\overline y$ coincide on $B$.  Hence  from inequality (\ref{ineqtree2}), $z_n(\mu)\leq \overline y(\mu)$ for $\mu \in B$.  We also have $Supp(\overline y)\cap B=Supp(\overline y_n)$. If  $K$ denote  this set, we  have $\overline z_n(\mu)=0$ if $\mu\in B\setminus K$. Now, since $K$ is finite and $\overline y$ takes only non-negative integer values, there are infinitely many $\overline z_n$'s which coincide on $K$. These $\overline z_n$'s coincide on $B$. \endproof

Let  $\overline z_{n_0}, \dots, \overline z_{n_{k}}, \dots$ be an infinite subsequence of $\overline z_n$'s such that the $\overline z_{n_{k}}$'s coincide on $B$. For each $k<\omega$ set $x_k:=\overline z_{n_0}\wedge   \overline z_{n_k}$. 

{\bf Claim 4.} For $k\geq 1$,  $x_k$ is the restriction of $z_{n_{0}}$ to $(\leftarrow \lambda_{n_{k}}[$. 

{\bf Proof of Claim 4.}

First, from our construction $z_{n_0}$ and $z_{n_k}$ coincide on $B$.  Next,  from Claim 2, $z_{n_0}$ and $z_{n_{k}}$ are $0$ on $(\leftarrow \lambda_{n_0}[\cap A$ and  $(\leftarrow \lambda_{n_k}[\cap A$ respectively. Since $n_0<n_k$, we have $\lambda _{n_k}< \lambda _{n_0}$.  It  follows that $z_{n_0}$ and $z_{n_k}$ coincide on $(\leftarrow \lambda_{n_k} [ $ and $z_{n_0}(\lambda_{n_k})=0 \not = z_{n_k}(\lambda_{n_k})$. The result follows. \endproof

From Claim 3, we immediately have:

{\bf Claim 4.} The sequence $x_0, \dots,x_k, \dots$ is strictly decreasing. 

 Consequently,  $T(X)$ contains an infinite strictly decreasing sequence. With that, the proof of the lemma is complete. 
\end{proof}
 
\begin{theorem} \label{thm:firstext}Every  countable ultrametric space $\M$ extends to an ultrametric commutative monoid  $\M^*$ having the same spectrum.  Moreover, if $\M$ has a well founded nerve,  $\M^*$ too. \end{theorem}
\begin{proof}
Let $\M$ be countable. Let $\lambda:= Spec(\M)\setminus \{0\}$ dually ordered. Then $\M$ isometrically embeds into $\omega^{[\lambda ]}$. Let $X$ be its image. Set $\M^{*}:= X^{*}$. \end{proof}

We may note that if $T$ has at least two elements then $M^{*}$ has infinitely many sons.

 \subsection{The path extension of an ultrametric space} We define the path extension $\mathbb Path(\M)$ of an ultrametric space $\M$. Its  elements are finite unions of chains in $(Nerv(\M) \supseteq)$. 
 \begin{notations}
Let $\M$  be an ultrametric space. Let  $r\in \R_{+}$, set $Nerv_{r}(\M):= \{B\in Nerv(\M): \delta(B)=r\}$. Let $\beta<\alpha \in \R^*_{+}\cup \{+\infty\}$.  Set $Nerv_{[\beta,\alpha[ }(\M):= \bigcup \{Nerv_{r}(\M): \beta\leq r<\alpha\}$ and set $Nerv_{< \alpha }(\M):=Nerv_{[0,\alpha[ }(\M)$. For $B\in Nerv_{< \alpha }(\M)$ set $$]_{\alpha}\leftarrow B]:= \{C\in Nerv_{< \alpha }(\M): C\supseteq B\}$$ 
\end{notations}
\begin{definition} A subset $\mathcal B$ of $Nerv(\M)$ is \emph{slim} if 
\begin{equation} \delta(B)=\delta (B') \Rightarrow B=B'
\end{equation}
for all $B, B'\in \mathcal B$
\end{definition}
Let $\mathcal B \subseteq Nerv(\M)$. We set $Spec (\mathcal B):=\{\delta (B): B\in \mathcal B\}$. If $Spec(\mathcal B)$ has a least element (w.r.t the order on the reals),  we denote it $\delta (\mathcal B)$. If moreover, $\mathcal B$ is slim, we denote by  $end(\mathcal B)$ the unique $B\in \mathcal B$ such that $\delta (B)= \delta(\mathcal B)$.  
Let $\mathcal B$ be a finite non-empty slim subset of $Nerv_{< \alpha }(\M)$.  Let $n:=\vert \mathcal B\vert $. We denote by  $\overline \beta:=(\beta_i)_{i<n}$  the unique enumeration of $Spec(\mathcal B)$ into a decreasing order 
$\beta_0>\dots >\beta_{n-1}$. The \emph{enumeration of} $\mathcal B$ is the sequence  $\overline{ \mathcal B} := (B_i)_{i<n}$  of elements of $\mathcal B$ such that $\delta(B_i)=\beta_i$.  We set $\beta_{-1}:= \alpha$ and we set:

$$]_{\alpha}\leftarrow \mathcal B]:=\bigcup  \{]_{\beta_{i-1}}\leftarrow B_i] :i<n\}$$.

\begin{definition} An \emph{$\alpha$-path}  is any subset $\mathcal I$ of $Nerv_{< \alpha }(\M)$ of the form $\mathcal I= ]_{\alpha}\leftarrow \mathcal B]$.  We denote by $L_{\alpha}$ the set of $\alpha$-paths. A finite set $\mathcal B$ \emph{generates} an $\alpha$-path $\mathcal I$ if $\mathcal I=  ]_{\alpha}\leftarrow \mathcal B]$.
\end{definition} 
\begin{fact}\label{subslim} 
If $\mathcal B$ is slim, $\mathcal B'\subseteq \mathcal B$ and $Spec(\mathcal B)=Spec(\mathcal B')$ then $\mathcal B=\mathcal B'$. 
\end{fact}
\begin{fact} \label{fact0} Let $B, B'\in Nerv_{< \alpha }(\M)$ such that $B'\supseteq B$. 
Then:  $$]_{\alpha}\leftarrow B]= ]_{\alpha}\leftarrow B'] \cup  ]_{\delta(B')}\leftarrow B]$$
\end{fact}
 \noindent{\bf Proof of Fact \ref{fact0}. } Observe that two balls containing $B$ are comparable w.r.t. inclusion. \endproof

\begin{fact} \label{fact1}Every  $\alpha$-path is slim.
\end{fact}

\begin{fact} \label{fact4}If $\mathcal B$   generates the  $\alpha$-path $\mathcal I$, then every finite  subset $\mathcal B'$ of $\mathcal I$ which contains $\mathcal B$ generates $\mathcal I$. 
  \begin{fact}\label{fact2}
 If $\mathcal B$ is a finite non-empty slim subset of $Nerv_{< \alpha }(\M)$ then $end(\mathcal B)= end(]_{\alpha}\leftarrow \mathcal B])$. \end{fact}

  \end{fact} \begin{fact} \label{fact3}A set $\mathcal B$ generates an $\alpha$-path $\mathcal I$ if and only if it satisfies the following conditions:
 \begin{enumerate}
\item $\mathcal B$ is a finite subset of $\mathcal I$.
 \item $end(\mathcal B)=end({\mathcal I})$.
 \item If $\overline {\mathcal B}$ is the enumeration of $\mathcal B$ then  for every $B\in \mathcal I$, if $\mathcal \beta_{i-1}>\delta(B)\geq \delta (B'_{i})$ for some $i<n$  then $B\supseteq B'_{i}$.
  \end{enumerate}
  \end{fact}
  
  \noindent{\bf Proof of Fact \ref {fact3}.} The three conditions stated are obviously necessary. Suppose that they hold. According to Fact \ref{fact1},  $\mathcal J:=]_{\alpha}\leftarrow \mathcal B]$ is well-defined. Condition (3) yields that $\mathcal I\subseteq \mathcal J$. For the converse, let  $i<n$ and $\mathcal A_i:=]_{\beta_{i-1}}\leftarrow B_i]$.  We prove that $\mathcal A_i \subset \mathcal I$. For that, let $\mathcal B'$ such that $\mathcal I= ]_{\alpha}\leftarrow \mathcal B']$ and  let $\overline {\mathcal B'}:= (B'_i)_{i<n'}$ be the corresponding enumeration of $\mathcal B'$.  Since $B_i\in \mathcal I$ there is some  $B'_{j}$  such that $\delta(B'_{j-1})>\delta(B_i)$ and $B_i\supseteq B'_j$. It follows that $]_{\delta(B'_{j-1})}\leftarrow B_i]\subseteq I$. If $\delta (B'_{j-1})\geq \delta(B_{i-1})$ we are done. If not, let $B'\in \mathcal A_i$ and $B'_{k}$  such that $\delta(B'_{k-1})>\delta(B')\geq  B'_k$. If $k=j$, $B'\supseteq B'_j$,  hence $B'\in \mathcal I$. If $k<j$ then  since $\mathcal B$ satisfies  Condition (3),  we have $B'_{k-1} \supseteq B_i$. From Fact \ref{fact0} this yields $B'\geq  B'_k$, hence $B'\in  \mathcal I$.  \endproof

  \begin{fact} \label{fact5}Let $\mathcal I$ be an $\alpha$-path. Then $\mathcal I\cap  Nerv_{<\beta}(\M)$ is a $\beta$-path provided that it is non empty and $\alpha >\beta$.
  
   \end{fact}\label{fact6} 
  \begin{fact} Let $\mathcal I$ be an $\alpha$-path. Then $\mathcal I\setminus Nerv_{<\beta}(\M)$ is an $\alpha$-path whenever  $\beta \in Spec (\mathcal I)$.
   \end{fact}  
   \noindent{\bf Proof of Fact \ref{fact6}. } Let $B\in \mathcal I$ such that $\delta(B)= \beta$. Let $\mathcal B$ which generates $\mathcal I$. According to Fact \ref{fact4} we may suppose that $B\in \mathcal B$. Use the definition of $ ]_{\alpha}\leftarrow \mathcal B]$ to conclude.\endproof

  \begin{fact}  \label{fact7} Let $\mathcal I$ be an $\alpha$-path and  $\mathcal J$ be a $\delta(\mathcal I)$-path. Then $\mathcal I \cup \mathcal J$  is an  $\alpha$-path.
\end{fact}

\begin{definition} A finite slim subset $\mathcal B$ of $Nerv_{<\beta}(\M)$ is \emph{pure} if two consecutive terms in the  enumeration of $\mathcal B$ are incomparable w.r.t inclusion. 
\end{definition}

\begin{lemma} Every $\alpha$-path $\mathcal I$ is generated by a unique pure set.
\end{lemma}
\begin{proof}
Let $\mathcal B$ a generating subset of $\mathcal I$ with minimum size and let $\overline {\mathcal B}:= (B_i)_{i<n}$ be its enumeration. Then $\mathcal B$ is pure. Indeed, suppose that $B_{i} $ and $B_{i+1}$ are comparable. Then $B_{i}\supset B_{i+1}$.  It follows that  $\mathcal B\setminus \{B_i\}$ satisfies the conditions of Fact \ref{fact3}. Hence it generates $\mathcal I$. This contradicts the minimality of the size of  $\mathcal B$. We show the uniqueness of $\mathcal B$ by induction on $n:= \vert \mathcal B\vert $. 

\begin{claim}
 Let $\mathcal I_0$ be the subset of $\mathcal I$ made of the elements $B$ such that:
\begin{enumerate}[{(a)}]
\item $]_{\alpha}\leftarrow B]\subseteq \mathcal I$
\item for every $B'\in \mathcal I$, if $\alpha >\delta(B')\geq \delta (B)$  then $B'\supseteq B$.
\end{enumerate}
Then 
$\mathcal I_0:= ]_{\alpha}\leftarrow B_0]$
\end{claim}
The proof is immediate and we omit it. 
Now, set $\mathcal I':= \mathcal {I}\setminus \mathcal {I}_0$. If $\mathcal {I'}=\emptyset$ we are done. If not, then $\mathcal {I'}$ is a $\beta_0$-path (Fact \ref{fact5}) and $\mathcal {I'}= ]_{\beta_0}\leftarrow \mathcal {B'}]$  where $\mathcal {B'}:= \mathcal B\setminus \{B_0\}$. Clearly $\mathcal {B'}$ has minimum size. Hence induction applies. The result follows. 
\end{proof}
\begin{notat}Let $\mathcal I$ be a slim set, let $a$ with $end (\mathcal I)< a$. Let $P_a(\mathcal I)$ be the set of $x\in Spec(\mathcal I)$ such that:

\begin{equation}
 x\leq \delta( \mathcal B)\leq \delta (\mathcal B')<a\Longrightarrow B\subseteq B'
 \end{equation}
 for all $B,B'\in \mathcal I$.
 If this set has a least element, we denote it by $\mu_a(\mathcal I)$.
 
 For an example, if $\mathcal I$ is  an $\alpha$-path, $\mu_{\alpha} (\mathcal I)=Max(Spec(\mathcal B))$ where $\mathcal B$ is the pure set generating $\mathcal I$. In this case, we set $l(\mathcal I):= \vert \mathcal B \vert$, $\mu (\mathcal I):= \mu_{\alpha} (\mathcal I) $, $init (\mathcal I):= B_0$ such that  $B_0\in \mathcal I$ and $\delta(B_0)= \mu(I)$. 
\end{notat}

Let $\mathcal I',\mathcal I''$ be two $\alpha$-paths. We set $\mathcal I'\leq_{\alpha} \mathcal I''$ if there is some $\beta\in Spec(\mathcal I'')$ such that $\mathcal I'= \mathcal I''\setminus \M_{<\beta}$. 
Let $ \mathcal B'$ and $\mathcal B''$ be the pure generating subsets of $\mathcal I'$ and $\mathcal I''$ respectively and let $\overline{\mathcal{B'}}:= (B'_i)_{i<n'}$ and $\overline{\mathcal{B''}}:= (B''_i)_{i<n''}$ be the corresponding sequences. Set $\overline{\mathcal{B'}}_{*}=(B'_i)_{i<n'-1}$ and $\overline{\mathcal{B''}}_{*}=(B''_i)_{i<n''-1}$.

Set $\overline{\mathcal{B'}}\leq_{\alpha} \overline{\mathcal{B''}}$ if $\overline{\mathcal{B'}}_{*}$ is a prefix of $\overline{\mathcal{B''}}_{*}$ and $B'_{n-1}\supseteq B''_{n-1}$.

\begin{fact} \label{fact9}
We have $\overline{\mathcal{B'}}\leq_{\alpha} \overline{\mathcal{B''}}$ if and only if $\mathcal I'\leq_{\alpha}\mathcal I''$.
\end{fact}

\begin{fact}\label {fact8}
$\mathcal I'\leq_{\alpha} \mathcal I''$ if and only if there is some $\delta(\mathcal I')$-path $\mathcal J''$ such that $I''= \mathcal I' \cup \mathcal J''$ \end{fact}
\begin{remark}
If $\mathcal I'\leq_{\alpha}\mathcal I''$ then $\mathcal I' \subseteq \mathcal I''$. But the converse does not necessarily holds. 
\end{remark}

For $\mathcal J\in  L_{\alpha}$, set $( \leftarrow \mathcal{J}]_{L_{\alpha}} := \{\mathcal I\in  L_{\alpha}:  \mathcal I\leq_{\alpha}\mathcal I\}$. 

\begin {fact} \label{fact9}Let $\mathcal {I'}, \mathcal {I''}\in ( \leftarrow \mathcal{J}]_{L_{\alpha}}$. Then: $\mathcal {I'}\leq_{\alpha} \mathcal {I}''$ if and only if $\delta (\mathcal I')\geq \delta(\mathcal {I''})$.  \end{fact}

\begin{fact} \label {fact10}The relation $\leq_{\alpha} $ is an order on the set $L_{\alpha}$ of $\alpha$-paths. For every $\mathcal J\in  L_{\alpha}$, the set  $( \leftarrow \mathcal{J}]_{L_{\alpha}}$ is linearly ordered. \end{fact}

\begin{fact} \label {fact11} If $(Nerv_{<\alpha}(\M), \supseteq)$ is well-founded, then $(L_{\alpha}, \leq_{\alpha})$ too. 
\end{fact}
 \noindent{\bf Proof of Fact \ref{fact11}. }
 Let $\mathcal{J}\in Nerv_{<\alpha}(\M)$. Observe that $Spec(\mathcal J)$ is a finite union of dually well ordered chains and apply Fact \ref{fact9}. \endproof

Let $\bot_{\alpha}$ be a set not belonging to $L_{\alpha}$. Extend the order $\leq_{\alpha}$ to   $\overline{L}_{\alpha}:= L_{\alpha}\cup \{\bot_{\alpha}\}$,  with the requirement that $\bot_{\alpha}\leq \mathcal I$ for all $\mathcal I \in L_{\alpha}$.
\begin{lemma}
Two elements $\mathcal I', \mathcal I''$ of $\overline{L}_{\alpha}$ have an infimum in  $\overline{L}_{\alpha}$ that we will denote $\mathcal I'\wedge_{\alpha} \mathcal I''$. 
\end{lemma}
\begin{proof}
If $\mathcal I'$ and $\mathcal I''$ are comparable, we have $\mathcal I'\wedge_{\alpha}\mathcal I''= Min\{\mathcal I',  \mathcal I''\}$. Otherwise, proceed by induction on $n:= l(\mathcal I')+l(\mathcal I'')$. Set 
$B'_0:= init(\mathcal I')$, $B''_0:= init(\mathcal I'')$, $\beta:= \delta(B'_0\cup B''_0)$, $B:= B'[B'_0\cup B''_0, \beta)$  and $\mathcal A:= ]_{\alpha}\leftarrow B'_0]\cap  ]_{\alpha}\leftarrow B''_0]$.  Hence $\mathcal A=]_{\alpha}\leftarrow B]$ if $\alpha>\beta$. 

{\bf Case 1.} $B'_0=B''_0$. In this case $B=B_0$. Set $\mathcal I'_{1}:= \mathcal {I' } \setminus \mathcal A$ and $\mathcal I''_{1}:= \mathcal {I'' } \setminus \mathcal A$. Since $\mathcal I'$ and $\mathcal I''$ are incomparable, $\mathcal I'_{1}$ and $\mathcal I''_{1}$ are non-empty. Hence $\mathcal I'_{1}, \mathcal I''_{1}\in  {L}_{\beta}$ and $l(\mathcal I'_{1}) +l( \mathcal I''_{1})=n-2$. Hence induction applies. Let $I_1:= \mathcal I'_{1}\wedge_{\beta}\mathcal I''_{1}$ in $\overline {L}_{\beta}$. If $I_1=\bot_{\beta}$ then $\mathcal I'\wedge_{\alpha}  \mathcal I''= \mathcal A$. If $I_1 \not=\bot_{\beta}$ then 
$\mathcal I'\wedge_{\alpha} \mathcal I''= \mathcal A\cup I_1$.

{\bf Case 2.} $B'_0\not =B''_0$. In this case, $\mathcal I'\wedge \mathcal I''= \mathcal A$ if $\alpha>\beta$. Otherwise $\mathcal I'\wedge_{\alpha} \mathcal I''= \bot_{\alpha}$.
 
\end{proof}

\begin{notations}Let $\alpha\in \R_{+}\cup \{+\infty\}$. Set $\delta(\bot_{\alpha})=\mu(\bot_{\alpha}):=\alpha$. For $ \mathcal I', \mathcal I''\in \mathbb \overline{L}_{\alpha}$, set $d_{\alpha}(\mathcal {I'}, \mathcal {I''}):= \delta(\mathcal I'\wedge_{\alpha} \mathcal I'')$. Let $\beta<\alpha$, an \emph{$(\alpha,\beta)$-path} is any $\alpha$-path $\mathcal I$ such that $\delta (\mathcal I)= \beta$. We denote $L_{\alpha, \beta}$ the set of $(\alpha, \beta)$-paths. We set $\overline {L}_{\alpha, \beta}:= L_{\alpha, \beta}\cup \{\bot_{\alpha}\}$. \end{notations}
\begin{lemma} \label{lem:distance}Let  $ \mathcal I', \mathcal I'', \mathcal J\in \mathbb \overline{L}_{\alpha}$. Then:
\begin{enumerate}
\item  $d_{\alpha}(\mathcal {I'}, \mathcal {I''})\leq \mu (\mathcal I'\wedge_{\alpha} \mathcal I'')$.
\item $d_{\alpha}(\mathcal {I'}, \mathcal {I''})\leq Max \{d_{\alpha}(\mathcal {I'}, \mathcal {J}), d_{\alpha}(\mathcal {I''}, \mathcal {J})\}$.\\
Moreover, if $ \mathcal I', \mathcal I'', \mathcal J\in \mathbb \overline{L}_{\alpha, \beta}$ then:
\item $d_{\alpha}(\mathcal {I'}, \mathcal {I''})=\beta$ if and only if $ \mathcal I'=\mathcal I''$.
\end{enumerate}
 \end{lemma}
 \begin{proof}
 \noindent  Item (1). According to our definitions of $\delta$ and $\mu$, we have  $\delta(\mathcal J)\leq  \mu(\mathcal J)$ for all $\mathcal J\in L_{\alpha}$. 
 
 \noindent Item (2). We have 
$\mathcal {I'}\wedge_{\alpha} \mathcal {J}\leq_{\alpha} \mathcal I'$ and $\mathcal {I''}\wedge_{\alpha} \mathcal {J}\leq_{\alpha} \mathcal I''$. Hence $\mathcal {I'}\wedge_{\alpha}\mathcal {I''}\wedge \mathcal J\leq_{\alpha} \mathcal {I'}\wedge_{\alpha} \mathcal I''$ (*). We have $\mathcal {I'}\wedge_{\alpha} \mathcal {J}\leq_{\alpha} \mathcal J$ and $\mathcal {I''}\wedge_{\alpha} \mathcal {J}\leq_{\alpha} \mathcal J$.  Hence, $\mathcal {I'}\wedge_{\alpha} \mathcal {J}$ and $\mathcal {I''}\wedge_{\alpha} \mathcal {J}$ are comparable (Fact \ref{fact10}). Suppose   $\mathcal {I'}\wedge_{\alpha} \mathcal {J}\leq_{\alpha} \mathcal {I''}\wedge_{\alpha} \mathcal {J}$. In this case, (*) yields  $\mathcal {I'} \wedge_{\alpha}\mathcal J \leq_{\alpha} \mathcal {I'}
\wedge_{\alpha} \mathcal I''$. Hence, with Fact \ref{fact9},  $d_{\alpha}(\mathcal {I'}, \mathcal {J})\geq d_{\alpha}(\mathcal {I'}, \mathcal {I''})$.

\noindent Item (3). 
Apply Fact \ref{fact9}.
 \end{proof}
 
 \begin{notat}
We set  $\M_{\alpha}:=(M, d\wedge \alpha)$ , where $d \wedge \alpha $ is defined by $d\wedge \alpha(x,y) := Min( \{d(x,y), \alpha\})$. For $x\in M$,  we set $\varphi_{\alpha}(x):= ]_{\alpha} \leftarrow \{x\}]$. We denote by $d$ be the restriction of $d_{\alpha}$ to $L_{\alpha, 0}$, we set    $Path_{\alpha} (\M):=L_{\alpha, 0}$ and $\mathbb Path_{\alpha}(\M):= (Path_{\alpha}(\M), d)$.
 \end{notat}
 
 \begin{theorem}\label{thm:pathext}Let $\M$ be an ultrametric space and $\alpha\in \R^*_{+}\cup \{+\infty\}$. Then:
 \begin{enumerate} 
\item $\mathbb Path_{\alpha}(\M)$ is an ultrametric space.
\item The map $\varphi_{\alpha}$ is an isometric embedding of  $\M_{\alpha}$ into $\mathbb Path_{\alpha}(\M)$. Moreover $Spec(\M_{\alpha}, x)=Spec (\mathbb Path_{\alpha} (\M), \varphi_{\alpha}(x))$ for every $x\in \M$.
\item $Spec(\M_{\alpha})=Spec(\mathbb Path_{\alpha}(\M))$.
 \end{enumerate}
 \end{theorem}
 \begin{proof}
\noindent Item (1). Apply item (2) and Item (3) of Lemma \ref{lem:distance}.

\noindent Item (2).

Let $x,y\in M$. Set $r:= d_{\alpha}(x,y)= Min \{d(x,y), \alpha\}$,  $X:=  \varphi_{\alpha}(x)\wedge_{\alpha}\varphi_{\alpha}(y)$. According to Lemma \ref{lem:distance},  we have $d(\varphi_{\alpha}(x), \varphi_{\alpha}(y))=\delta(X)$. Let $B$ be the least ball in $Nerv(\M)$ containing $x$ and $y$. If $d(x,y)\geq \alpha$, that is $r = \alpha$, then  $X:=\bot_{\alpha}$, hence $\delta(X)=\alpha$. If not,   $X=]_{\alpha} \leftarrow B]$,  in which case $\delta(X)=r$.  In both cases  $d(\varphi_{\alpha}(x), \varphi_{\alpha}(y))=r$. This proves that $\varphi_{\alpha}$ is an embedding. 
Let $\mathcal I\in Path_{\alpha}(\M)$ such that $d(\varphi_{\alpha}(x), \mathcal I)= r$. Hence for $Y:= \mathcal{I}\wedge_{\alpha} \varphi_{\alpha}(x)$ we have $\delta (Y)=r$.  Pick $y\in \M$ such that $\{y\}\in \mathcal I$. We have $d(x,y)\geq  \delta (Y)$. Thus $d(x,y)\geq \alpha$ if $\delta(Y)=\alpha$. If not,  there is $B\in Nerv(M)$ such that $Y=]_{\alpha}\leftarrow B]$. In this case $\delta(B)=r$. Thus  $r\in Spec(\M_{\alpha}, x)$ .    

\noindent  Item(3). $\delta(\mathcal I)\in Spec (\M_{\alpha})$ for every $\mathcal I\in L_{{\alpha}, 0}$ hence $Spec(\mathbb Path_{\alpha}(\M))\subseteq Spec_{\alpha}(\M)$. The reverse inclusion  follows from Item 2. \end{proof}

 \begin{notat} We denote by $\mathbb Path(\M)$ the space  $\mathbb Path_{+\infty}(\M)$. We call it the \emph{path extension} of $\M$.  For  $r\in \mathbb R^*_{+}$  and $\mathcal I\in L_{+\infty, r}$, we set $ \mathcal I \ast L_{r, 0}:= \{Iæ\cup J: J\in L_{r, 0}\}$.
 \end{notat}
 \begin{lemma} \label {lem:ball}Let $B\subseteq Path(\M)$. Then $B$ is a non trivial member of $Nerv(\mathbb Path(\M))$ if and only if there are $r\in \mathbb R_{+}\setminus \{0\}$ and $\mathcal I\in L_{+\infty, r}$ such that $B= \mathcal I \ast L_{r, 0}$.
 \end{lemma}
 \begin{proof}Let $\mathbb P:= \mathbb Path(\M)$. Suppose that $B=\mathcal I \ast L_{r, 0}$. Let $X\in B$ and $Y\in Path(\M)$ such that $d(X,Y)\leq r$. We have  $\delta(X\wedge_{+\infty} Y)\leq r$. Thus 
 $\mathcal I\leq_{+\infty} Y$, or  equivalently  $Y= X\cup Y'$ for some $Y'\in L_{r, 0}$, that is  $Y\in B$. Hence $B\in  Nerv_r(\mathbb P)$. Conversely, let $B\in  Nerv_r(\mathbb P)$. Then there are $X, Y\in B$ such that $d(X,Y)=r$. Set $\mathcal I:= X\wedge_{+\infty} Y$ and $B':= \mathcal I \ast L_{r, 0}$. Since, as we have just seen, $B'$ is a ball or radius $r$, and since $B'$ contains $B$, we have $B=B'$.  \end{proof}
 \begin{corollary} Two members of $Nerv(\mathbb Path(\M))$ are isometric if and  only if they have the same diameter.
 \end{corollary}
\begin{proof} 
Let $B, B'\in Nerv_r(\mathbb Path(\M))$. We may suppose $r\not = 0$. From Lemma \ref{lem:ball}, $B= \mathcal I \ast L_{r, 0}$ and $B'= \mathcal I' \ast L_{r, 0}$. For $\mathcal J\in B$ set $f(\mathcal J):= \mathcal I'\cup( \mathcal J\setminus \mathcal I)$. Then $f$ defines an isometry of $B$ onto $B'$.
\end{proof}
\begin{lemma}
Let $r\in Spec(\M)\setminus \{0\}$. If some $B\in Nerv_r(\M)$ has infinitely many sons, then every $B'\in Nerv_r(\mathbb Path(\M))$ has infinitely many sons.
\end{lemma} 
\begin{proof}Let $B'\in Nerv_r(\mathbb Path(\M))$. According to Lemma \ref{lem:ball},   $B= \mathcal I \ast L_{r, 0}$ for some $\mathcal I\in L_{+\infty, r}$. For $x\in M$, set $\varphi_r(x):=  ]_{r} \leftarrow \{x\}]$ and $\theta(x):= \mathcal I \cup \varphi_r(x)$. Then, as  one can readily see, $\theta$ is an isometry from $B$ into $B'$. Since these two balls have the same diameter, $B'$ has as many sons as $B$. 
\end{proof}

 \subsection{Extension of a compatible operation to the path extension of an ultrametric space}\label{section:extension}
 In this section we extend a compatible operation on an ultrametric space $\M$ to its path extension $\mathbb Path(\M)$. The path we follow  is motivated by the following observation. If the operation on $\M$ is a kind of addition, then viewing members of $\mathbb Path(\M)$  as kind of  piecewise linear maps each defined on a subdivision of an interval $[\beta, \alpha[$, the natural idea to add two maps, $f$ and $g$  is to take a common refinement, and add the maps on the intervals of the refinement. But, as in our case,  it is possible that one map, say $f$,  is undefined on some interval $I$ of the refinement, we are forced to look at the values of $f$ outside $I$, and this makes the definition of sum a bit more complicated. 

Our result is this:

 \begin{theorem}\label{thm:compat}Let $\M$ be an ultrametric space and $\alpha\in \R^*_{+}\cup \{+\infty\}$. Suppose that there is a compatible  binary operation $+$ on $M$ . Then
there is a compatible operation $+_{\alpha}$ on  $\mathbb Path_{\alpha}(\M)$ such that:
\begin{equation}
\varphi_{\alpha}(x+y)= \varphi_{\alpha}(x)+_{\alpha}\varphi_{\alpha}(y)
\end{equation}
 for every $x, y\in M$.
  If moreover $+$ is associative, resp. commutative, resp. has a neutral element  then $+_{\alpha}$ too. And if $0$ is the neutral element for $+$ then $\varphi_{\alpha}(0)$ is the neutral element for $+_{\alpha}$.  
 \end{theorem}

The proof will occupy the rest of this section. 

We  extend successively  the operation $+$ to $Nerv_{<\alpha}(\M)$, to $\check{ S}_{\alpha, \beta}$ and to $Path_{\alpha, \beta}(\M)$.

The extension to $Nerv_{<\alpha}(\M)$ is immediate. 
\begin{lemma} Let $\M$ be an ultrametric space  endowed with a compatible binary operation $+$. Then for all $x,x',y,y'\in M$:
 \begin{equation} \label{ineqmetr2}
d(x+x',y+y')\leq Max \{d(x,y), d(x', y')\}
\end{equation}
\end{lemma}
Inequality\ref{ineqmetr2} asserts that  $+$ is a non-expansive map from $\M\times \M$ equipped with the $\ell^{\infty}$ metric.  
Let $B, B'\subseteq M$. Set $B\check{+}B':= \{x+x': x\in B, x'\in B'\}$.
\begin{fact} If $\M$ is an ultrametric space, and $B,B'$ are non-empty then:
\begin{equation}
\delta(B\check{+}B')=Max\{\delta(B), \delta(B')\}
\end{equation}
Moreover, $\delta(B\check{+}B')$ is attained whenever $\delta(B)$ and $\delta(B')$ are attained. 
\end{fact}
\begin{notat}
If  $B,B'$ are two bounded subsets of $M$, we denote by $B + B'$ the least member of $Nerv(\M)$ containing $B\check{+}B'$. 
\end{notat}

Ordered by reverse of the inclusion, $Nerv(\M)$ is a meet-lattice. In lattice  terms,  $B + B'= \bigwedge \{\{x+x'\}: x\in B , x'\in B'\}$. This extend to $Nerv_{<\alpha}(\M)$ provided that a least element $\bot _{\alpha} $ is added.

The following relationship between the operation $+$, the meet and $\delta$ is the clue for a proof of  the theorem.  
\begin{lemma} \label{heyting} Let $\M$ be an ultrametric space  endowed with a compatible binary operation $+$ and $\alpha \in \R^*_{+}\cup \{+\infty\}$. Then 
\begin{enumerate}
\item $C+(B\wedge_{\alpha} B')=(C+B)\wedge_{\alpha} (C+B')$

\item $(B\wedge_{\alpha} B')+C=(B+C)\wedge_{\alpha} (B+C)$

\item $\delta (B+B')=Max\{\delta(B), \delta(B')\}$
\end{enumerate}

for all $B,B'\in Nerv_{<\alpha} (\M)$.

\end{lemma}
One may note that,  conversely, an operation on $Nerv_{<\alpha}(\M)$ satisfying the three conditions of the lemma come from a compatible operation on $\M$.
 To go further we  need the following:
  
 \begin{notations} Let $\beta\in \R_{+}$, $\alpha \in \R_{+}\cup \{+\infty\}$ such that $\beta<\alpha>\beta$. Let $ S_{\alpha, \beta}$ be the set of slim subsets $\mathcal B$ of $Nerv_{[\beta, \alpha[}(\M)$  such that $\beta\in Spec (\mathcal B)$ and let $\check{S}_{\alpha, \beta}$ the subset of those which are finite. 

For a finite subset $X$ of $[\beta, \alpha [$ containing $\beta$ and  $a \in  [\beta,  \alpha [$,  set $X(a):= Max \{x: \beta \leq x\leq a\}$. For $\mathcal B\in \check{  {S}}_{ \alpha, \beta}$,  let $\mathcal B(a)$ be the unique $B\in \mathcal B$ such that $\delta(B):= Spec(\mathcal B)(a)$.   
\end{notations}
 
 \begin{lemma}\label{lem:sumslim}
 Let $n$ be a positive integer. Let  $\oplus: Nerv_{[\beta, \alpha[}(\M)^n\rightarrow Nerv_{[\beta, \alpha[}(\M)$ be an $n$-ary operation on $Nerv_{[\beta, \alpha[}(\M)$.  For every $(\mathcal B_i)_{i<n}\in \check{{S}}^{n}_{\alpha, \beta}$, set  
\begin{equation}\label{aq:sum}
\overline{ \oplus}_{i<n} \mathcal {B}_{i}:= \{\oplus_{i<n}\mathcal {B}_{i}(a): a\in [\beta,  \alpha[\}.
\end{equation}

\noindent Suppose that:
 \begin{equation}\label {eq:ballcompat}
  \delta(\oplus_{i<n} B_i)= Max(\{\delta(B_i):i<n\})
 \end{equation}
for all $(B_i)_{i<n}\in Nerv_{[\beta, \alpha[}(\M)^n$. 

\noindent (1) Let   $(\mathcal B_i)_{i<n}\in \check{ {S}}^{n}_{\alpha, \beta}$
then: 
\begin{enumerate} [{(a)}] 

\item  $Spec(\overline \oplus_{i<n} \mathcal {B}_{i})= \cup_{i<n} Spec (\mathcal B_{i})$
\item $\overline \oplus_{i<n} \mathcal {B}_{i}\in\check{ S}_{\alpha, \beta}$. 
\item $(\overline \oplus_{i<n} \mathcal {B}_{i})(a)= \oplus_{i<n} \mathcal {B}_{i}(a)$ for all $a\in[\beta,  \alpha[$.
\end{enumerate}

\noindent  Suppose moreover that :
  \begin{equation}\label {eq:ballcompat2}
  \oplus_{i<n} B_i\subseteq \oplus_{i<n}B'_i
 \end{equation}
   for all $(B_i)_{i<n}, (B'_i)_{i<n}\in Nerv_{[\beta, \alpha[}(\M)^n$ such that $B_i\subseteq B'_i$ for every $i<n$.

\noindent (2) Let  $(\mathcal J_{i})_{i<n}\in {L}^{n}_{\alpha, \beta}$ and  $(\mathcal B_i)_{i<n}, (\mathcal B'_i)_{i<n}\in \check{ {S}}^{n}_{\alpha, \beta}$.
 
If $\mathcal B_i$ and  $\mathcal B'_i$ generates $\mathcal J_i$ for every  $i<n$ then $\overline \oplus_{i<n} \mathcal {B}_{i}$ and $\overline \oplus_{i<n} \mathcal {B'}_{i}$ generate the same member of ${{L}}_{\alpha, \beta}$. 

\end{lemma}
\begin{proof}
Let  $a\in [\beta,  \alpha[$, From inequality (\ref {eq:ballcompat}) and the definition of $\mathcal {B}_{i}(a)$ we have
\begin{equation}\label{eq:trivial}
\delta( \oplus \mathcal {B}_{i}(a))=Max\{\delta(\mathcal {B}_{i}(a)): i<n\} \leq a
\end{equation} 
(1) Item  (a) follows immediately. 
%
Item (b). Suppose that $\delta( \oplus \mathcal {B}_{i}(a))=\delta(\oplus \mathcal {B}_{i}(a'))$. Suppose that $a\not =a'$. W.l.o.g. we may suppose $a<a'$. Let $i<n$. From inequality (\ref{eq:trivial})
we have $\delta(\mathcal {B}_{i}(a))\leq \delta( \oplus \mathcal {B}_{i}(a))\leq a$, hence  $\delta(\mathcal {B}_{i}(a))=\delta(\mathcal {B}_{i}(a'))$. Since $\mathcal B_i$ is slim, this yields $\mathcal {B}_{i}(a)=\mathcal {B}_{i}(a')$. From which we get $\oplus \mathcal {B}_{i}(a)=\oplus \mathcal {B}_{i}(a')$, proving that $\overline{ \oplus}_{i<n} \mathcal {B}_{i}$ is slim. 
 
 Item (c). From inequality  (\ref{eq:trivial}) we have $\delta(\oplus _{i<n}\mathcal {B}_{i}(a))=a$ if and only if $a\in  \cup_{i<n} Spec (\mathcal B_{i})$. Since $\overline{ \oplus}_{i<n} \mathcal {B}_{i}$ is slim, this yields $(\overline \oplus_{i<n} \mathcal {B}_{i})(a)= \oplus_{i<n} \mathcal {B}_{i}(a)$.

\noindent (2)To simplify notations, set $\mathcal B:= \overline {\oplus}_{i<n} \mathcal {B}_{i}$, $A:= Spec(\mathcal B)$, $A_i:=Spec(\mathcal B_i)$, $\mathcal J:= ]_{\alpha}\leftarrow \mathcal B]$, and define similarly $\mathcal B',A',A'_i, \mathcal J'$. 
We prove our affirmation in two steps.
%


\begin{claim} \label{claim:notobvious} If $\mathcal B_i\subseteq \mathcal B'_i$ for all $i<n$. then $\mathcal J'=\mathcal J$.
\end{claim}
\noindent {\bf Proof of Claim \ref{claim:notobvious}.}

\noindent Subclaim 1. Let $a\in [\beta, \alpha[$.  Then $\mathcal B(a)\subseteq \mathcal B'(a)$. 

\noindent Proof. Since $\mathcal B_i\subseteq \mathcal B'_i$ we have $A_i\subseteq A'_i$, hence we  have $A_{i}(a)\leq A'_{i}(a)\leq a$. That is 
$\delta(\mathcal B_{i}(a))\leq  \delta (\mathcal B'_{i}(a))\leq a$. Since $\mathcal B_{i}$ generates $\mathcal J_{i}$, $]_{a}\leftarrow \mathcal B_i(a)]\subseteq \mathcal J_i$.   Item 3 of Fact \ref{fact3} yields that 
$\mathcal B _{i}(a) \subseteq \mathcal B'_{i}(a)$. According to inequality \ref {eq:ballcompat2},  $\mathcal B(a)\subseteq \mathcal B'(a)$.\endproof

\noindent Subclaim 2. $\mathcal B'\subseteq J$. 

\noindent Proof. Since  $\mathcal B$ generates $\mathcal J$, $]_{a}\leftarrow \mathcal B(a)]\subseteq \mathcal J$. Since $\mathcal B(a)\subseteq \mathcal B'(a)$ we have $\mathcal B'(a)\subseteq ]_{a}\leftarrow \mathcal B(a)]$ and since $\mathcal B$ generates $\mathcal J$, $\mathcal B'(a)\subseteq ]_{a}\leftarrow \mathcal B(a)]$\endproof.

\noindent Subclaim 3. $\mathcal B\subseteq \mathcal B'$.  

\noindent Proof. Let $a\in A$. We have $\delta (\mathcal B(a))=a$ and $\delta (\mathcal B'(a))\leq a$.  From Subclaim 1, we have $\mathcal B(a)\subseteq \mathcal B'(a)$. Hence $\mathcal B(a)=\mathcal B'(a)$, proving that $B(a)\in \mathcal B'$. \endproof

Since $\mathcal B\subseteq \mathcal B'\subseteq J$ and $\mathcal B$ generates $\mathcal J$, it follows from Fact  \ref{fact4} that $\mathcal B'$ generates $\mathcal J$. Since $\mathcal B'$ generates $\mathcal J'$, we have $\mathcal J'=\mathcal J$ as claimed. \endproof

 Let  $(\mathcal {B''})_{i<n}$ be  a family of finite slim sets, each $\mathcal B''_i$ generating $\mathcal J_i$. Let   $\mathcal J''$ be the corresponding $(\alpha, \beta)$-path.
  \begin{claim}\label{claim:obvious}
 $\mathcal J=æ\mathcal J'' $ 
\end{claim}

\noindent {\bf Proof of Claim \ref{claim:obvious}.}
Let  $(\mathcal {B'})_{i<n}$, where $\mathcal B'_i:= \mathcal B_i\cup \mathcal B''_i$ for all $i<n$. According to Claim \ref{claim:notobvious} we have $\mathcal J= \mathcal J'$ and $\mathcal J=\mathcal J''$. The result follows. 
\endproof

This completes the proof of the lemma.
\end{proof}

\begin{definition}Let  $\oplus: Nerv_{[\beta, alpha[}(\M)^n\rightarrow Nerv_{[\beta, alpha[}(\M)$ be an $n$-ary operation  satisfying conditions (\ref{eq:ballcompat}) and (\ref{eq:ballcompat2}). Let $(\mathcal J_i)_{i<n}\in L^n_{\alpha, \beta}$. We set $\overline \oplus (\mathcal J_i)_{i<n}:= ]_{\alpha} \leftarrow \overline \oplus (\mathcal B_i)_{i<n}]$ where $(\mathcal B_{i})_{i<n}\in \check{S}^n_{\alpha, \beta}$ is such that  $\mathcal J_i= ]_{\alpha} \leftarrow\mathcal B_i]$. According to Lemma \ref{lem:sumslim}, this does not depends upon the choice of $(\mathcal B_{i})_{i<n}$. \end{definition}

\begin{corollary}\label{cor: phi} Let $+$ be a compatible operation on $\M$.  The  operation + defined on $\mathbb Path_{\alpha}(\M)$ satifies
\begin{equation}
\varphi_{\alpha}(x+y)= \varphi_{\alpha}(x)+_{\alpha}\varphi_{\alpha}(y)
\end{equation}
\end{corollary}
\begin{proof}
$\varphi_{\alpha}(x)= ]_{\alpha} \leftarrow \{x\}]$ thus the pure slim set $\mathcal B$ generating $\varphi_{\alpha}(x)$ reduces to a singleton (namely $\{\{x\}\}$). The result follows.\end{proof}

\begin{corollary}\label{cor:assoc}
Let $+$ be a binary operation on $Nerv_{[\beta, \alpha[}(\M)$ satisfying conditions (\ref{eq:ballcompat}) and (\ref{eq:ballcompat2}). If this operation is associative, then the extensions of this operation to   $ \check{ S}_{\alpha, \beta}$ and to $L_{\alpha, \beta}$  are associative.  
\end{corollary}
\begin{proof}
We  set $\oplus_{i<2}B_i:=B_0+B_1$ and $\oplus_{i<3}B_i:= B_0+B_1+B_3$. Applying the definitions given in Lemma \ref{lem:sumslim}, we have with obvious notations:
\begin{equation}\label {claim:assoc1}
(\mathcal {B}_{1}+\mathcal {B}_{2})+\mathcal {B}_{3}=\mathcal {B}_{1}+\mathcal {B}_{2}+\mathcal {B}_{3}
\end{equation}
The associativity of the extension of $+$ to $ \check{ S}_{\alpha, \beta}$ follows. 
The definitions of the corresponding operations on $L_{\alpha, \beta}$ yield the same formula.\end{proof}

\begin{lemma}\label{lem: extpath} Let $\M$ be a metric space and  $+$ be a  compatible binary operation on $\M$. Then its extension  to $L_{\alpha, \beta}$ is compatible:

\begin{equation}\label {eq:comp}
d_{\alpha}(\mathcal J+\mathcal {I'}, \mathcal J+\mathcal {I''})= d_{\alpha}(\mathcal {I'}, \mathcal {I''})=d_{\alpha}(\mathcal {I'}+\mathcal J, \mathcal I'+\mathcal {J})\end{equation}
for all $\mathcal {I'}, \mathcal {I''},  \mathcal J$ in $L_{\alpha, \beta}$.
\end{lemma}
\begin{proof}
Let $\gamma:= d_{\alpha}(\mathcal {I'}, \mathcal {I''})$ and $\gamma':= d_{\alpha}(\mathcal J+\mathcal {I'}, \mathcal J+\mathcal {I''})$. 
Let $\mathcal {B'}, \mathcal {B''}, \mathcal C$ be pure sets  belonging to  $\check{S}_{\alpha, \beta}$ and generating respectively $\mathcal {I'},\mathcal {I''}$ and  $\mathcal {J}$. Since $\mathcal {B'}$ and  $\mathcal {B''}$ are pure,   $\mathcal {B'}(a)= \mathcal {B''}(a)$ for all $a \in [\gamma,  \alpha[$. Thus $\mathcal {C}(a)+\mathcal {B'}(a)=\mathcal {C}(a)+\mathcal {B''}(a)$ for all $a \in [\gamma,  \alpha[$.  Since $\mathcal {J+I'}$ and $\mathcal {J+I''}$ are respectively generated by $\mathcal {C+B'}$ and $\mathcal {C+B''}$, this yields  $\gamma'\leq \gamma$. 

For the converse, suppose $\gamma'<\gamma$. Set $b:=\mu_{\gamma}(\mathcal B)$, 
$b':=\mu_{\gamma}(\mathcal B')$, $c:=\mu_{\gamma}( \mathcal C)$ and $d:=\mu_{\gamma}( \mathcal D)$, $d':=\mu_{\gamma}( \mathcal D')$ where $\mathcal D, \mathcal D'$ are the two pure sets generating $\mathcal J+\mathcal I$ and $\mathcal J+\mathcal I'$ respectively.  Set $e:=Max\{b,b',c,c',d,d', \gamma'\}$.
Since $\gamma', d,d'\leq e$, we have 
\begin{equation}
(\mathcal C(e)+\mathcal B(e))\wedge_{\alpha}(\mathcal C(e)+\mathcal B(e))\in \mathcal J+\mathcal I'\wedge_{\alpha} \mathcal J+\mathcal I''
\end{equation}
Hence $\delta(\mathcal C(e)+\mathcal B(e))\wedge_{\alpha}(\mathcal C(e)+\mathcal B(e))\leq e$.

On an other hand, since $Max \{b,b'\}\leq e<\gamma$, we have $\delta(\mathcal B(e)\wedge_{\alpha}\mathcal B(e))=\gamma$.
Thus $\delta (\mathcal C(e)+\mathcal B(e)\wedge_{\alpha}\mathcal B(e))= Max\{\delta (\mathcal C(e)),\delta(\mathcal B(e)\wedge_{\alpha}\mathcal B(e))\}=\gamma$.
 From the distributivity property stated in Lemma \ref{heyting}, we have 
\begin{equation}\label {eq:distrib}
(\mathcal C(e)+\mathcal B(e))\wedge_{\alpha}(\mathcal C(e)+\mathcal B(e))= \mathcal C(e)+\mathcal B(e)\wedge_{\alpha}\mathcal B(e)
\end{equation}
This yields a contradiction.\end{proof}

\noindent {\bf Proof of Theorem \ref{thm:compat}.} Let $+$ be a compatible operation on $\M$. According to Lemma \ref{heyting} it extends to an operation on $Nerv(\M)$ satisfying conditions \ref{eq:ballcompat} and \ref{eq:ballcompat} of Lemma \ref{lem:sumslim}.  Hence it extends to an operation on $\mathbb Path_{\alpha}(\M)$, which according to Lemma \ref{lem: extpath} is compatible. According to Corollary \ref {cor:assoc}, this operation is associative provided that the original one is associative. 

\endproof

\subsection{Indivisibility properties of the path extension of an ultrametric space.}

\begin{lemma}\label {eqcompatible2} If a metric space $\M$ can be endowed with a compatible binary operation $+$ then
\begin{enumerate}
\item For every $a\in M$, $\M$ embeds into $\M_{\restriction M(a)}$.
\item If $M$ has infinitely many sons then for each son $B$ of $M$ there are infinitely many sons $B'$ such that $\M_{\restriction B}$ embeds  into $\M_{\restriction B'}$.
\end{enumerate}
\end{lemma} 
 \begin{proof}

\begin{claim} \label{claim:compat}
\begin{equation}\label{eqcompatible3} Spec (\M, a) \subseteq Spec (\M, a+b)
\end{equation}
for all $a, b \in \M$
\end{claim}
{\bf Proof of Claim \ref{claim:compat}.} 
Let $r\in Spec(\M, a)$. Let $x\in M$
 such that $d(a, x)=r$. Then  from  (\ref{eqcompatible1}), $d(a+b, x+b)=r$ proving that $r\in Spec(\M, a+b)$.
\endproof

Item 1.  Let $a \in M$. Let $T_a: M\rightarrow M$ defined by setting $T_a(x):= a+x$. Since $+$ is compatible, $T$ is an embedding of $\M$. From Claim \ref{claim:compat}, we have $T_a(x)\in M(a)$ for every $x\in M$, as required. 

Item 2 Let $B'$ be a son of $M$ and $b\in B'$, select $a_B''$ in each son $B''$ distinct from $B'$. Then the images of $B'$ by the $T_{a_B''}$'s are in different sons. 
\end{proof}

\begin{theorem} \label{doubleextension1} 
The path extension $\mathbb Path(\M)$ of a countable ultrametric monoid  $\M$  is  spec-endogeneous and hereditary indivisible provided that $(Nerv(\M), \supseteq )$ is well founded and has infinitely many sons \end{theorem}
\begin{proof}
It  suffices to prove  that $\mathbb Path(\M)$ satisfies conditions (1), (3)  and (4) of Theorem \ref{main}.

Since  $\mathbb Path (\M)$ is an extension of $\M$,  it has infinitely many sons. Hence, from Lemma \ref {eqcompatible2}: 
\begin{enumerate}[{(a)}]

 \item $\mathbb Path( \M)$ embeds into $\mathbb Path( \M)_{\restriction Path( \M)(x)}$ for every $x\in Path(\M)$.

 \item For each  son $B$ of $Path(\M)$ there are infinitely many sons $B'$ such that $\mathbb Path(\M)_{\restriction B}$ embeds  into $\mathbb Path(\M)_{\restriction B'}$.
 
From Theorem \ref{thm:compat} these properties hold for every $B\in Nerv(\mathbb Path (\M))$ replacing $\mathbb Path (\M)$. 
\end{enumerate}

In particular, every non trivial member of $Nerv(\mathbb Path(\M))$ has infinitely many sons, this is condition (3). Since two members of $Nerv(\mathbb Path(\M))$ with the same diameter are isometric, property (ii)(a) of Proposition \ref{prop:updirected} holds.  Since property (b) above is property (ii)(b) of Proposition \ref{prop:updirected}, we get that (i)(a) and (i) (b) holds, that is conditions (1) and (4) hold.    \end{proof}


\begin{thebibliography}{999}


\bibitem{DLPS} C.Delhomm\'e, C.Laflamme, M.Pouzet, N.Sauer, Divisibility of countable metric spaces, preprint, july 2005, 33pp., http://arxiv.org/math.CO/051025, to appear in {\em Europ. J. of Combinatorics} (2006). 
\bibitem {DLPSHom} C.Delhomm\'e, C.Laflamme, M.Pouzet, N.Sauer, A note on homogeneous ultrametric spaces,  february  2007. 
 \bibitem{El}   M. El-Zahar, N.W. Sauer,    A Game for Vertex
Partitions,  {\em  Discrete Math.,} {\bf 291}(2005), 99-113.
\bibitem {fraisse} Fra\"\i ss\'e,  
Sur l'extension aux relations de quelques propri\'et\'es des ordres. {\em Annales scientifiques de l'\'Ecole Normale Sup\'erieure},  S\'er. 3, 71 no. 4 (1954), p. 363-388.   
\bibitem{Fra}   R. Fra\"\i ss\'e,   Theory of Relations, Revised Edition,
in: {\em  Studies in Logic and the Foundations of Mathematics},
 {\bf 145},  North Holland 2000
\bibitem{hjorth} G. Hjorth, An oscillation theorem for groups of isometries, {\em manuscript} Dec. 31, 2004, 28 pp.
%
\bibitem{KPT} A.S. Kechris, V.G. Pestov and S. Todorcevic, Fra\"\i ss\'e limits, Ramsey theory, and topological
dynamics of automorphism groups. {\em Geom. Funct. Anal.},  15:106--189, 2005.  
\bibitem {lemin2}A.Lemin, The category of ultrametric spaces is isomorphic to the
category of complete, atomic, tree-like, and real graduated
lattices LAT*, {\em  Algebra universalis},  {\bf 50} (2003),  35?49.
\bibitem{nesetril} J.Nesetril , Metric spaces are Ramsey. {\em European Journal of
Combinatorics}, 28 (2007) 457-468.

\bibitem{vanthe}
L.Nguyen Van Th\'e.
\newblock  Big Ramsey degrees and indivisibility  in classes of ultrametric spaces, preprint July 5, 2005, pp13.
\bibitem {pestov} V. Pestov, Ramsey--Milman phenomenon, Urysohn metric spaces, and extremely amenable
groups. {\em Israel Journal of Mathematics} 127:317-358, 2002. Corrigendum, ibid. 145:375-379,
2005.
\bibitem{pouzet} M.Pouzet, Relation impartible, {\em Dissertationnes}, {\bf 103} (1981), 1-48.

\bibitem{sauer} N.W. Sauer,    Canonical vertex partitions, 
{\em Combinatorics Probability and Computing.} {\bf 12}(2003) 671-704.


\bibitem{urysohn} P.Urysohn, Sur un espace m\'etrique universel, {\em Bull. Sci. Math.} (2) {\bf 51} (1927) 43-64.

\end{thebibliography}
\end{document}